\magnification=\magstep1
\def\refer#1{ \medskip  \parindent=34pt \hang \noindent
\rlap{[#1]\hfil} \hskip 30pt }

\baselineskip=14pt

\font\tenmsb=msbm10
\font\sevenmsb=msbm7
\font\fivemsb=msbm5
\newfam\msbfam
\textfont\msbfam=\tenmsb
\scriptfont\msbfam=\sevenmsb
\scriptscriptfont\msbfam=\fivemsb
\def\Bbb#1{{\fam\msbfam\relax#1}}

\def\de{{\delta}}
\def\ff{{\infty}}
\def\R{{\Bbb R}}
\def\E{{{\Bbb E}\,}}
\def\P{{\Bbb P}}

\def\Z{{\Bbb Z}}
\def\F{{\cal F}}

\def\lam{{\lambda}}

\def\al{{\alpha}}
\def\grad{{\nabla}}
\def\proof{{\medskip\noindent {\bf Proof. }}}
\def\longproof#1{{\medskip\noindent {\bf Proof #1.}}}
\def\qed{{\hfill $\square$ \bigskip}}
\def\subsec#1{{\bigskip\noindent \bf{#1}.}}

\def\cite#1{{[{#1}]}}

\def\norm#1{\Vert #1 \Vert}

 \def\qq {\qquad}
\def\frac#1#2{{#1\over #2}}

\def\wt{\widetilde}
\def\ol{\overline}

\def\ni{\noindent }
\def\ms{\medskip}
\def\bs{\bigskip}
\def\cl#1{\centerline{#1}}

\parindent=30pt

\def\square{{\vcenter{\vbox{\hrule height.3pt
           \hbox{\vrule width.3pt height5pt \kern5pt
              \vrule width.3pt}
           \hrule height.3pt}}}}

\def\tfrac#1#2{{\textstyle {#1\over #2}}}

\def\choose#1#2{{\left(\matrix{#1\cr#2\cr}\right)}}
\def\tlint{{- \kern-0.85em \int \kern-0.2em}}  
\def\dlint{{- \kern-1.05em \int \kern-0.4em}}  

\def\Bas{1}
\def\BKa{2}
\def\BKb{3}
\def\BKc{4}
\def\BR{5}
\def\BP{6}
\def\Da{7}
\def\Db{8}
\def\LGa{9}
\def\LGb{10}
\def\LGc{11}
\def\La{12}
\def\Lw{13}
\def\MR{14}
\def\Pt{15}
\def\Ra{16}
\def\Rb{17}
\def\Rc{18}
\def\Sp{19}
\def\Va{20}

\cl{\bf AN ALMOST SURE INVARIANCE PRINCIPLE}
\ms
\cl{\bf FOR RENORMALIZED INTERSECTION LOCAL TIMES}
\bigskip

\vskip0.4truein
\centerline{{\bf Richard F.~Bass}
\footnote{$\empty^1$}{\rm Research partially supported by NSF grant
DMS-0244737.}
\quad  and \quad {\bf Jay Rosen}
\footnote{$\empty^2$}{\rm Research partially supported by grants  from
the NSF and from PSC-CUNY.} }

\vskip1truein

{\bf Abstract}. Let $\wt \beta_k(n)$ be the number of self-intersections of
order $k$, appropriately renormalized, for a mean zero  random walk
$X_n$ in $\Z^2$ with $2+\delta $ moments. On  a suitable  probability
space we can construct $X_n$ and a planar Brownian motion $W_t$ such
that for each $k\geq 2$
$$ |\wt \beta_k(n)-\wt \gamma_k(n)|=O(n^{ -a}), \qq \hbox{a.s.}$$
for some $a>0$ where
$\wt
\gamma_k(n)$ is the renormalized self-intersection local time of order
$k$ at time 1 for the Brownian motion $W_{nt}/\sqrt n$.

\bs
\ni Subject classifications: Primary 60F17; Secondary 60J55

\bs
\ni Key words: intersection local time, invariance principle, random walks,
strong invariance, self-intersection local time

\bigskip
\ni{Short title:  } Self-intersection local time

\vfill\eject

\subsec{1. Introduction}

If $\{W_{ t}\,;\,t\geq 0\}$ is a planar Brownian motion with density
$p_{ t}( x)$,  set
$\gamma_{1,\epsilon}( t)=t$ and for $k\geq 2$ and $x=(x_{ 2},\ldots,x_{ k}
)\in (\R^{2})^{k-1}$ let
$$
\gamma_{k,\epsilon}(t,x)=\int_{ 0\leq t_{ 1}\leq\cdots\leq t_{ k}<t}
\prod_{ i=2}^{ k} p_{ \epsilon}(W_{ t_{ i}}-W_{ t_{ i-1}}-x_{ i})dt_{ 1}\cdots
dt_{ k}.
$$ When $x_{ i}\neq 0$ for all $i$ the limit
$$
\gamma_{k}(t,x)=\lim_{ \epsilon\rightarrow 0}\,\gamma_{k,\epsilon}(t,x)
$$ exists and for any bounded continuous function $F( x)$ on $\R^{2( k-1)
}$ we have
$$  \int   F( x)\gamma_{k}(t,x)\,dx=\int_{ 0\leq t_{ 1}\leq\cdots\leq t_{
k}<t} F(W_{ t_{ 2}}-W_{ t_{ 1}},\ldots,W_{ t_{ k}}-W_{ t_{ k-1}})dt_{
1}\cdots dt_{ k}.  \eqno(1.1)  $$ (Here we may arbitrarily specify that
$\gamma_{k}(t,x)=\infty$ if any $x_{ i}= 0$.) When
$x_{ i}\neq 0$ for all
$i$ define the renormalized intersection local  times as
$$
\wt \gamma_{k}(t,x)=\sum_{A\subseteq \{ 2,\ldots,k\}}( -1)^{|A|}
\Big(\prod_{ i\in A} \frac{1}{\pi}\log ( 1/|x_{i}|)\Big)\gamma_{k-|A|}(t,x_{
A^{ c}})
$$ where
   $x_{A^c}=(x_{i_1},
\ldots, x_{i_{k-|A|}})$ with
$i_1<i_2<\cdots <i_{k-|A|}$ and $i_j\in \{2, \ldots, k\}-A$ for each
$j$, that is, the vector $(x_2, \ldots, x_k)$ with all terms that 
have indices in
$A$ deleted. It is known that the $\wt \gamma_{k}(t,x)$ have a continuous
extension to all $\R^{ 1}_{ +}\times \R^{ k-1}$; see
\cite{\BKb}.

Renormalized self-intersection local time was originally  studied by
Varadhan
\cite{\Va} for its role in quantum field theory. In  Rosen
\cite{\Rc} it is shown that
$\wt\gamma_{k}( t,0)$ can be characterized as the  continuous process of
zero quadratic variation in the decomposition of a natural Dirichlet process.
Renormalized intersection local time turns out to be the right tool for the
solution of certain ``classical'' problems such as the asymptotic expansion
of the area of the Wiener sausage in the plane and the range of
   random walks, \cite{\BR}, \cite{\LGa}, \cite{\LGb}. For further work on
renormalized self-intersection local times see Dynkin
\cite{\Da}, Le Gall \cite{\LGc}, Bass and Khoshnevisan
\cite{\BKb}, Rosen
\cite{\Rb} and Marcus and Rosen \cite{\MR}.

Let $\xi_i$ be i.i.d. random variables with values in
$\Z^2$ that are mean 0, with covariance matrix equal to the identity, and
with $2+\delta$ moments.  Let us suppose the
$\xi_i$ are symmetric and are strongly aperiodic.  Let $X_n$ be the random
walk, that is, $X_n=\sum_{i=1}^n \xi_i$. Let $p(n,x,y)$ be the transition
probabilities.   Let $B_1(n,x)=n$ and for $x\in \Z^2$ set
$$B_2(n,x)=
\sum_{0\leq i_{ 1}<i_{ 2}\leq n}
   1_{(X_{ i_{ 2}}=X_{ i_{ 1}}+x)}.$$ More generally, for
 $x=(x_2,\ldots, x_{ k})\in (\Z^2)^{ k-1}$ let \
$$B_k(n,x) =\sum_{0\leq i_{ 1}<i_{ 2}<\ldots<i_{ k}\leq n}
\prod_{j=2}^{k} 1_{(X_{ i_{ j}}=X_{ i_{ j-1}}+x_j)}.$$ Note that
$B_k(n,x)=0$ for all $n<k-1.$

With $e_1=( 1,0)$, let
$$ G(x)=\sum_{n=1}^\infty [p(n,0,x)-p(n,0,e_1)],$$
  and set
   $G_{ n}(x)=G(x)-G(\sqrt n e_1)$. Let
$$\wt B_k(n,x)=\sum_{A\subset\{2, \ldots,k\}} (-1)^{|A|}
\Big(\prod_{i\in A} G_{ n}(x_i)\Big) B_{k-|A|}(n,x_{A^c}).
\eqno (1.2)
$$ In particular we have
$$\wt B_2(n,x)=B_2(n,x)- G_{ n}(x)n.
\eqno (1.3)
$$

Finally we define the  renormalized intersection local  times for our random
walk by
$$\wt \beta_k(n,x)=\frac{1}{n} \wt B_k(n,x\sqrt n).\eqno (1.4)$$ In
particular we have
$$\wt \beta_2(n,x)=\frac{1}{n}B_2(n,x\sqrt n)- G_{ n}(x\sqrt n).
\eqno (1.5)
$$ We note from P12.3 of \cite{\Sp} that for $x\neq 0$
$$  \lim_{ n\rightarrow \infty}G_{ n}(x\sqrt n)=
\lim_{ n\rightarrow \infty}G(x\sqrt n)-G(\sqrt n e_1) =\frac{1}{\pi}\log (
1/|x|).  \eqno (1.6)
$$

We know we can find a version of our random walk and a Brownian motion
$W_t$ such that
$$\sup_{s\leq 1} |X^n_s-W^n_s|=o(n^{-\zeta}), \qq a.s. \eqno (1.7)$$  for
some $\zeta>0$ where
$X^n_t= X_{[nt]}/\sqrt n$ and $W^n_t= W_{nt}/\sqrt n$; see,
\cite{\BP}, Theorem 3,  for example. Let $\gamma_k(1,x,n)$ and
$\wt
\gamma_k(1,x,n)$ be the intersection local times and renormalized
intersection local  times up to time 1  of order $k$, resp.,  for the Brownian
motion
$W^n_t$. In this paper we prove the following theorem.

\proclaim Theorem 1.1. Let $X_{ n}=\xi_{ 1}+\cdots+\xi_{ n}$ be a random
walk in
$\Z^2$, where the $\xi_i$ are i.i.d., mean 0, with covariance matrix equal to
the identity, with
$2+\delta$ moments for some $\delta>0$, symmetric, and strongly
aperiodic. On a suitable probability space we can construct $\{X_{
n}\,;\,n\geq 1\}$ and a planar Brownian motion
$\{W_{ t}\,;\,t\geq 0\}$ and we can find $\eta>0$ such that for each
$k\geq 2$
$$|\wt \beta_k(n,0)-\wt \gamma_k(1,0,n)|=o(n^{-\eta}), \qq a.s.$$

For related work see \cite{\BKc}, \cite{\BR}, \cite{\Ra}.

We give a brief overview of the proof. There is an equation similar to (1.1)
when
$\gamma_k$ is replaced by $\wt \gamma_k$, and also when it is replaced
by $\wt \beta_k$. Since by (1.7) we have $X^n_s$ close to
$W^n_s$ for
$n$ large, we are able to conclude that
$\int F(x) \wt \gamma_k(1,x)\,dx$ is close to $\int F(x) \wt
\beta_k(n,x) \, dx$ for $n$ large. If $F$ is smooth and has integral 1, then
by the continuity of $\wt \gamma_k(t,x)$ in $x$, which is proved in
\cite{\BKb}, we see that $\int F(x)\wt \gamma_k(1,x) \, dx$ is not far from
$\wt \gamma_k(1,0)$. If we had a similar result for $\wt
\beta_k$, we would then have that
$\int F(x) \wt \beta_k(n,x)\, dx$ is not far from $\wt \beta(n,0)$, and we
would have our proof. So our strategy is to obtain  good estimates on $|\wt
\beta_k(n,x)-\wt \beta_k(n,0)|$. Because of the rate of convergence in (1.7),
it turns out we are able to avoid having to find the sharpest estimates on the
difference, which simplifies the proof considerably.

Our main tool in obtaining the desired estimates is Proposition 3.2. This
proposition may be of independent interest. It has been known for a long
time that one way of proving $L^p$ estimates for a continuous increasing
process is to prove corresponding estimates for the potential. It is 
not as well
known that one can do the same for continuous processes of bounded
variation provided one has some control on the total variation; see, e.g.,
\cite{\Bas} or \cite{\BKb}.  Proposition 3.2 is the discrete time analogue of
this result, and is proved in a similar way. Unlike the continuous time
version, here it is also necessary to have control on the differences of
successive terms.

Section 2 has some estimates on the potential kernel for random walks in
the plane, while Section 3 has the proof of the stochastic calculus results we
need. Theorem 1.1 in the case when $k=2$ is proved in Section 4, with the
proofs of some lemmas postponed to Section 5. We treat  the case
$k=2$ separately for simplicity of exposition. The description of the
potentials of intersection local times of random walks in the $k>2$ case is a
bit different than in the $k=2$ case and this is described in Section 6.
Theorem 1.1 in the $k>2$ case is proved in Section 7, with the proofs of
some lemmas given in Sections 8 and 9. Finally in Section 10 we give an
extension of our results to $L^2$ convergence, and more importantly, make
a correction to the proof of one of  the propositions in
\cite{\BKb}. An Appendix contains the detailed proof of that correction. Throughout this paper we use
the letter
$c$ to denote finite positive constants whose exact value is unimportant and which may vary
from line to line.

\subsec{2. Estimates for random walks}

In this section we prove some estimates for the potential kernel of a random
variable. See the
  forthcoming book by Lawler \cite{\Lw} for further information. Let
$\ol G$ be the potential kernel for $X$. Recall that in 2 dimensions, since
$X$ is recurrent, the potential kernel is defined somewhat differently  than
in higher dimensions, and is defined by
$$\ol G(x)=\sum_{n=0}^\infty [p(n,0,x)-p(n,0,e_1)],$$ where
$e_1=(1,0)$. (Note $e_1$ can be replaced by any fixed point.) For us it will
be more  convenient to work with
$$ G(x)=\sum_{n=1}^\infty [p(n,0,x)-p(n,0,e_1)],$$ which, since
$p(0,0,0)=1$ and
$p(0,0,e_1)=0$, differs from $\ol G(x)$ by $1_{\{0\}}(x)$. By Spitzer
\cite{\Sp}, p.~75, we have
$$p(n,0,x)\leq c/n. \eqno (2.1)$$  By \cite{\BKc}, Proposition 2.1, if the
$\xi_i$ are strongly aperiodic, then
$$|p(n,0,x)-p(n,0,y)|\leq \frac{c|x-y|}{n^{ 3/2}}.
\eqno (2.2)$$

\proclaim Proposition 2.1. Suppose the $\xi_i$ have  $2+\delta$ moments.
Then $G(x)$ exists  and $|G(x)|\leq c(1+\log^+ |x|)$.

\proof Using (2.2), we have that
$$|G(0)|\leq \sum_{n=1}^\infty \frac{c}{n^{3/2}}$$ is finite.  The rest of
the assertions follow from
$$\eqalign{|G(x)|&\leq \sum_{n=1}^{|x|^2} \frac{c}{n}
+\sum_{n=|x|^2+1}^\infty
\frac{c|x-e_1| }{n^{3/2}}\cr &\leq c+ c\log |x| +
c\frac{|x|}{(|x|^2+1)^{1/2}}.\cr}$$
\qed

\proclaim Proposition 2.2. For some $c<\ff$
$$|G(x)-G(y)|\leq c\Big(\frac{|x-y|}{(1+|x|)\land (1+|y|)}\Big)^{2/3}, \qq
x,y\in \Z^2.$$

\proof By \cite{\Sp}, P7.10,
$$p(j,0,x)\leq \frac{c}{|x|^2}. \eqno (2.3)$$ Since $p(j,0,0)\leq 1$, then we
have
$$p(j,0,x)\leq \frac{c}{1+|x|^2}. \eqno (2.4)$$ Suppose $0<|x|\leq |y|$. Let
us set $R$ in a moment. Using (2.4) for $j\leq R$ and (2.2) for
$j>R$, we have that
$$|G(x)-G(y)|\leq \sum_{j=1}^R \frac{c}{1+|x|^2} +\sum_{j=R+1}^\infty
\frac{c|x-y|}{j^{3/2}}
\leq \frac{cR}{1+|x|^2} +\frac{c|x-y|}{R^{1/2}}.$$ If we select
$R$ so that
$$\frac{R}{1+|x|^2}  =\frac{|x-y|}{R^{1/2}},\qquad{\rm  i.e.,  }\qquad
R^{3/2}=(1+|x|^2)(|x-y|),$$ the result follows. Since
$G(0)$ is finite and
$|G(x)|\leq c\log(1+|x|)\leq |x|^{ 2/3}$, the result holds when either
$x$ or $y$ is 0, as well.
\qed

\proclaim Lemma 2.3.   For some constant $\kappa$ and any
$\rho<\de/2$,
$$G(x)=\kappa+\tfrac{1}{\pi} \log (1/|x|)+O( |x|^{-\rho}), \qq x\in \Z^2. $$

\proof Let us begin with the proof of Proposition 3.1 in \cite{\BKa}.  We
have for
$\delta>0$
$$\Bigl| e^{i\al\cdot x/\sqrt n}-1 -\frac{\al\cdot x}{\sqrt n} -\frac{|\al\cdot
x|^2}{2n}\Bigr|\leq c\Bigl|\frac{\al\cdot x}{\sqrt n}\Bigr|^ {2+\delta}.$$ So
if $\phi$ is the characteristic function of a random vector with finite
$2+\delta$ moments, mean 0, and the identity as its covariance matrix, then
$$\phi(\al/\sqrt n)=1-\frac{|\al|^2}{2n}+ E_1(\al,n),$$ with
$$|E_1(\al,n)|\leq c(|\al|/\sqrt n)^{2+\delta}. \eqno (2.5)$$ Applying this
also for  the characteristic function of a standard normal vector,
$$e^{-|\al|^2/2n}=1-\frac{|\al|^2}{2n}+E_2(\al,n),$$ where
$E_2(\al,n)$ has the same bound as $E_1(\al,n)$. If we use this in place of
the display in the middle of page 473 of \cite{\BKa}, we obtain
$$I_1^{(n)}\leq c n^{-\delta/2} (\log n)^{(4+\delta)/2}.$$ So if
$E(n,x)=|p(n,0,x)-(2\pi n)^{-1} e^{-|x|^2/2n}|$, following the proof  in
\cite{\BKa} we obtain
$$\sup_x E(n,x)\leq cn^{-1-(\delta/2)} (\log n)^{(4+\delta)/2}.$$ Let us
choose $\delta'<\delta$. We then have
$$\sup_x E(n,x)\leq cn^{-1-(\delta'/2)}.\eqno (2.6)$$

Recall $G(x)=\sum_{k=1}^\infty [p(k,0,x)-p(k,0,e_1)]$. It is shown in the
proof of Theorem 1.6.2 in \cite{\La} that for some constant $\kappa$
$$ \sum_{k=1}^\infty [q(k,0,x) -q(k,0,e_1)]=\kappa+\tfrac{1}{\pi}
\log(1/|x|) +O(|x|^{-1}),$$ where $q(k,x,y)=(2\pi k)^{-1} e^{-|x-y|^2/2k}$.
Thus, to prove the lemma, it suffices to prove
$$ \sum_{k=1}^\infty |p(k,x, 0) -q(k,x,0)|= O(|x|^{-\rho}) \eqno (2.7)$$ for
any
$\rho<\de/2$.

To establish (2.7), use \cite{\Pt}, p.~60 to observe that
$$p(k,x,0)\leq \P(|X_k|>|x|)\leq
\frac{\E |X_k|^{2+\delta}}{|x|^{2+\delta}}\leq c\frac{k^{1+\delta/2}}
{|x|^{2+\delta}}$$ and a similar estimate is easily seen to hold for
$q(k,x,0)$. Therefore, using (2.6) and setting $R= |x|$,
$$\eqalign{\sum_{k=1}^\infty |p(k,x, 0) -q(k,x,0)| &\leq
\sum_{k=1}^{R} [p(k,x,0) + q(k,x,0)]\cr &\qq +\sum_{k=R}^\infty  |p(k,x,
0) -q(k,x,0)|\cr &\leq \sum_{k=1}^{R}
c\frac{k^{1+\delta/2}}{|x|^{2+\delta}} +\sum_{k=R}^\infty
c\frac{1}{k^{1+\delta'/2}}\cr &\leq
c\frac{R^{2+\delta/2}}{|x|^{2+\delta}}+c \frac1{R^{\delta'/2}}\cr &\leq
c|x|^{-\rho}.\cr}$$
\qed

\subsec{3. Stochastic calculus}

   We will use the following propositions; these may be of independent
interest. Propositions 3.1 and 3.2 and their proofs are the discrete time
analogues of Propositions 6.1 and 6.2 of \cite{\BKb}.

\proclaim Proposition 3.1. Let $A_n$ be an adapted increasing sequence  of
random variables with $A_0=0$ and $A_\infty =\sup_n A_n$ finite.
Suppose that
$$Y=\sup_n (A_n-A_{n-1})$$ and
$W$ is a random variable such that
$$\E[A_\infty-A_n\mid \F_n]\leq \E[W\mid \F_n]$$ for all $n$. Then for
each integer
$p$  larger than 1 there  exists a constant $c$ such that
$$\E A_\infty^p\leq c p^p (\norm{W}_p+\norm{Y}_p)^p.$$

\proof Since $A_n$ is increasing,
$$pA_n^{p-1} \geq A_n^{p-1}+
A_n^{p-2}A_{n-1}+\cdots+A_nA_{n-1}^{p-2}+ A_{n-1}^{p-1}.$$ Multiplying
by $A_n-A_{n-1}$, we obtain
$$(A_n-A_{n-1})pA_n^{p-1}\geq A_n^p-A_{n-1}^p.$$ Summing over
$n$ we obtain
$$p\sum_{ n=1}^{ \infty} (A_n-A_{n-1}) A_{n}^{p-1}\geq
A_\infty^p.\eqno (3.1)$$ On the other hand, applying the general
summation formula
$$A_{  \infty}B_{  \infty}=
\sum_{ n=1}^{ \infty}A_{n}( B_{ n} -B_{n-1}) +\sum_{ n=1}^{ \infty}
(A_{n}-A_{n-1}) B_{ n-1}$$ with
$B_{ n}=A_{n+1}^{ p-1}$ we obtain
$$\eqalignno{\sum_{ n=1}^{ \infty}(A_{n}-A_{n-1})\, A_{n}^{ p-1}&= A_{
\infty}^{ p}-\sum_{ n=1}^{ \infty}A_{n} ( A_{n+1}^{ p-1}-A_{n}^{
p-1})&(3.2)\cr &=\sum_{ n=1}^{ \infty} (A_\infty-A_n)( A_{n+1}^{
p-1}-A_{n}^{ p-1}) +A_\infty A_{1}^{ p-1}.\cr}$$  Here we used the fact that
$$\sum_{ n=1}^{ \infty} A_\infty
(A_{n+1}^{p-1}-A_{n}^{p-1})=A_\infty^p-A_\infty A_{1}^{ p-1}.$$
Combining (3.1) and (3.2) we obtain
$$A_\infty^p\leq p\sum_{n=1}^\infty (A_\infty
-A_{n})(A_{n+1}^{p-1}-A_{n}^{p-1}) +pA_\infty A_{1}^{ p-1}.
\eqno (3.3)$$

Now suppose for the moment that $Y$ is bounded and
$A_{n}=A_{n_0}$  for all $n\geq n_0$ for some $n_0$. We have
$$\eqalignno{& \sum_{n=1}^\infty (A_\infty
-A_{n})(A_{n+1}^{p-1}-A_{n}^{p-1})& (3.4)\cr &=\sum_{n=1}^\infty
(A_\infty-A_{n+1})(A_{n+1}^{p-1}-A_{n}^{p-1}) +\sum_{n=1}^\infty
(A_{n+1} -A_{n}) (A_{n+1}^{p-1}-A_{n}^{p-1}).\cr}$$ and
$$\sum_{n=1}^\infty (A_{n+1} -A_{n}) (A_{n+1}^{p-1}-A_{n}^{p-1})\leq
Y   A_\infty^{p-1}.
$$ But, $A_{ 1}\leq A_\infty$ and also $A_{ 1}\leq Y$ so that
$$A_\infty A_{1}^{ p-1}\leq Y   A_\infty^{p-1}.      $$

  We write
$$\eqalign{\E \sum_{n=1}^\infty
(A_\infty-A_{n+1})(A_{n+1}^{p-1}-A_{n}^{p-1})& =\E
\sum_{n=1}^\infty
\E[A_\infty-A_{n+1}\mid \F_{n+1}](A_{n+1}^{p-1}-A_{n}^{p-1})\cr
&\leq
\E
\sum_{n=1}^\infty
\E[W\mid  \F_{n+1}] (A_{n+1}^{p-1}-A_{n}^{p-1})\cr &= \E
\sum_{n=1}^\infty W (A_{n+1}^{p-1}-A_{n}^{p-1})\cr &\leq  \E [W
A_\infty^{p-1}].\cr}$$ Therefore using H\"older's inequality,
$$\E A_\infty^p \leq p (\norm{W}_p+\norm{Y}_p) (\E
A_\infty^p)^{1-\frac1{p}}.$$ Our temporary assumptions on $A$ allow us
to divide both sides by $(\E A_\infty^p)^{1-\frac1{p}}$ to obtain our result
in this special case.

In general, look at
$$ A'_n=\sum_{ j=1}^{ n} ((A_j-A_{j-1})\land K)$$ and apply the above to
$A''_n=A'_{n\land n_0}$; note that $A''$ will satisfy the hypotheses with
the same $W$  and $Y$. Then let
$K\uparrow \infty$ and next
$n_0\uparrow \infty$ and use monotone  convergence.
\qed

\proclaim Proposition 3.2. Suppose $Q^1_n$ and $Q^2_n$ are two adapted
nonnegative increasing sequences. Suppose
$Q_n=Q^1_n-Q^2_n$,
$H_n+Q_n$ is a martingale that is 0 at time 0, $$Z=\sup_n |H_n|,$$
$$Y=\sup_n [(Q^1_n-Q^1_{n-1})+(Q^2_n-Q^2_{n-1})],$$ and
$$W=Q^1_\infty+Q^2_\infty.$$ Then there exists $c$ such that for $p>1$
$$\eqalignno{\E \sup_n |Q_n|^{2p}& \leq cp^{4p} \Big[ \E Z^{2p} +(\E
Z^{2p})^{1/2} (\E W^{2p})^{1/2}&(3.5)\cr &\qq\qq +(\E Y^{2p})^{1/2} (\E
W^{2p})^{1/2}\Big].\cr}$$

\proof There is nothing to prove unless $\E W^{2p}<\infty$. Since
$\sup_n Q_n\leq W$, all the random variables that follow will satisfy the
appropriate integrability conditions. Let us temporarily assume that there
exists $n_0$ such that
$Q_n^i=Q_{n_0}^i$ if $n\geq n_0, i=1,2$.

Let
$$V_m=\E[Q_\infty-Q_m\mid \F_m], \qq M_m=\E[Q_\infty\mid
\F_m].$$ Note that
$V_\infty=0$, $M_m$ is a martingale, and
$Q_m=M_m-V_m$. In fact, in view of our temporary assumption,
$V_m=0$ if $m\geq n_0$.

Our first observation is that since
$$V_m=\E[Q_\infty-Q_m\mid\F_m]=\E[H_m-H_\infty\mid
\F_m],$$ then
$$|V_m|\leq 2\E[Z\mid \F_m]. \eqno (3.6)$$ By Doob's inequality,
$$\E \sup_n V_n^p\leq c\E Z^p. \eqno (3.7)$$

We will use

\proclaim Lemma 3.3.
$$\E \Big(\sum_{n=0}^\infty (M_{n+1}-M_n)^2\Big)^p
\leq c(\E [ZW+YW]^p+ \E Z^{2p}).\eqno (3.8)$$

This lemma will be proved shortly. We first show how Proposition 3.2 follows
from this lemma. By the Burkholder-Davis-Gundy inequalities, we obtain
$$\E \sup_n |M_n|^{2p} \leq c p^{4p}(\E[ZW+YW]^p +\E Z^{2p}).
\eqno (3.9)$$ Combining with (3.7) and the fact that
$Q_m=M_m-V_m$ and then using Cauchy-Schwarz  completes the proof of
Proposition 3.2 in the special case where the $Q^i$ are constant from some
$n_0$ on. In the general case, let $\ol Q^i_n=Q^i_{n\land n_0}$ for
$i=1,2$, obtain (3.5) for $\ol Q_n=\ol Q^1_n-\ol Q^2_n$, let $n_0\to
\infty$, and apply monotone convergence.
\qed

{\bf Proof of Lemma 3.3.}   We now prove (3.8). Algebra shows that
$$V_\infty^2-V_m^2=\sum_{n=m}^\infty
(V_{n+1}-V_n)^2+2\sum_{n=m}^\infty V_n(V_{n+1}-V_n) .\eqno(3.10)$$
(Note that the sums are actually finite because $V_m=0$ if $m\geq n_0$.)
Recalling
$Q=M-V$ and $V_\infty=0$, we have
$$\eqalignno{\sum_{n=m}^\infty (M_{n+1}-M_n)^2&=\sum_{n=m}^\infty
(V_{n+1}-V_n)^2+\sum_{n=m}^\infty (Q_{n+1}-Q_n)^2\cr &\qq
+2\sum_{n=m}^\infty (Q_{n+1}-Q_n)(M_{n+1}-M_n)\cr
&=-V_m^2-2\sum_{n=m}^\infty V_n(V_{n+1}-V_n) +\sum_{n=m}^\infty
(Q_{n+1}-Q_n)^2\cr &\qq +2\sum_{n=m}^\infty
(Q_{n+1}-Q_n)(M_{n+1}-M_n)\cr &\leq -2\sum_{n=m}^\infty
V_n(V_{n+1}-V_n) +\sum_m (Q_{n+1}-Q_n)^2\cr &\qq
+2\sum_{n=m}^\infty (Q_{n+1}-Q_n)(M_{n+1}-M_n)\cr
&=:-2S_1+S_2+2S_3.&(3.11)\cr}$$

  We now take the  conditional expectation with respect to $\F_m$.
$$\eqalign{\E[S_1\mid \F_m]&=\E\Big[\sum_{n=m}^\infty
V_n(V_{n+1}-V_n)\mid
\F_m\Big]\cr &=\E\Big[\sum_{n=m}^\infty \E[V_{n+1}-V_n\mid
\F_n] V_n\mid
\F_m\Big]\cr &=-\E\Big[\sum_{n=m}^\infty \E[Q_{n+1}-Q_n\mid
\F_n] V_n\mid
\F_m\Big]\cr &=-\E\Big[\sum_{n=m}^\infty V_n (Q_{n+1}-Q_n)\mid
\F_m\Big]\cr &=-\E\Big[\sum_{n=m}^\infty V_{n+1} (Q_{n+1}-Q_n)\mid
\F_m\Big]\cr &\qq +\E\Big[\sum_{n=m}^\infty
(V_{n+1}-V_n)(Q_{n+1}-Q_n)
\mid \F_m\Big]\cr &=:I_1+I_2.\cr}$$

Since $V_{n+1}\leq 2\E[Z\mid F_{n+1}]$ by (3.6), we have
$$|I_1|\leq 2\E\Big[\sum_{n=m}^\infty Z|Q_{n+1}-Q_n|\mid
\F_m\Big]
\leq 2\E[ZW\mid \F_m].$$

Recalling that $V_n=M_{ n}-Q_n$ we see that
$$I_2=-\E[S_2\mid \F_m]+\E[S_3\mid \F_m].$$

   Since
$$\sum_{n=m}^\infty (Q_{n+1}-Q_n)^2\leq YW\eqno (3.12)$$ we have
$$\E[S_2\mid \F_m] \leq \E[YW\mid \F_m].$$ Let
$J=\E[\sum_{n=m}^\infty (M_{n+1}-M_n)^2\mid \F_m]$. By
Cauchy-Schwarz and (3.12),
$$|\E[S_3\mid \F_m]|\leq J^{1/2}\Big(\E[YW\mid
\F_m]\Big)^{1/2}.$$
   We therefore conclude
$$J\leq c\E[ZW+YW\mid \F_m]+ c_1 J^{1/2} \Big(\E[YW\mid
\F_m]\Big)^{1/2}.\eqno (3.13)$$  Using the inequality $A^{1/2}
x^{1/2}\leq (A+x)/2$ with $x=J$ and
$A=c_1^2 \E[YW\mid \F_m]$, we see that
$$J\leq c\E[ZW+YW\mid \F_m]+ c_1^2\Big(\E[YW\mid
\F_m]\Big)/2+ J/2 \eqno (3.14)$$ and therefore
$$J=\E\Big[\sum_{n=m}^\infty (M_{n+1}-M_n)^2\mid \F_m\Big]\leq
c\E[ZW+YW\mid \F_m]. \eqno (3.15)$$

We have $|Q_{n+1}-Q_n|\leq Y$ and so
$$\E\sup_n |Q_{n+1}-Q_n|^p\leq c\E Y^p.\eqno (3.16)$$ Using (3.7),
(3.16) and the fact that $Q_m=M_m-V_m$ and $Y\leq (YW)^{ 1/2}$ we
then have that
$$\E \sup_n |M_{n+1}-M_n|^{2p}\leq cp^{2p} (\E Z^{2p}+\E
(YW)^{p}).\eqno (3.17)$$ (3.8) then follows using
   (3.15) and Proposition 3.1 with $A_{ n}=\sum_{j=1}^n
(M_{j}-M_{j-1})^2$.\qed

\subsec{4. The $k=2$ case}

\proclaim Proposition 4.1. If
$$\wt U_2(n,x)=\sum_{i=0}^{n-1} G_n(X_n-X_i-x),
$$ then
$$M_n=\wt U_2(n,x)+\wt B_2(n,x)$$ is a martingale with $M_0=0$.

\proof If
$$U_2(n,x)=\sum_{i=0}^{n-1} G(X_n-X_i-x),
$$ we have $\wt U_2(n,x)=U_2(n,x)-nG( e_{1}\sqrt{n})$ so that
$$M_n= U_2(n,x)+ B_2(n,x)-nG(x).$$ Abbreviating
$\bar B_n= B_2(n,x)-nG(x)$ we have
$$\bar  B_n- \bar B_{n-1}=\sum_{i=0}^{n-1} 1_{(X_n=X_i+x)}-G(x).$$
  So 
$$\eqalign{\E[\bar   B_n -\bar  B_{n-1}+G(x)\mid \F_{n-1}]&
=\sum_{i=0}^{n-1}
\P(X_n-X_{n-1}+X_{n-1}-X_i=x\mid \F_{n-1})\cr &=\sum_{i=0}^{n-1}
p(1,0,X_{n-1}-X_i-x)\cr}\eqno (4.1)$$ Abbreviating
$U_n= U_2(n,x)$ we have
$$ U_n- U_{n-1}-G(x)=
\sum_{i=0}^{n-1} [G(X_n-X_i-x)-G(X_{n-1}-X_i-x)].$$ Now
  for any $i\leq n-1$
$$\eqalignno{\E[G(X_n-X_i-x)\mid
\F_{n-1}]&=\E[G(X_n-X_{n-1}+X_{n-1}-X_i-x)\mid
\F_{n-1}]&(4.2)\cr &=\sum_y G(y+X_{n-1}-X_i-x) \P(X_n-X_{n-1}=y
\mid
\F_{n-1})\cr &=\sum_y G(y+X_{n-1}-X_i-x) \P(X_n-X_{n-1}=y)\cr
&=\sum_y G(y+X_{n-1}-X_i-x)p(1,0,y)\cr &=P_1G(X_{n-1}-X_i-x)
\cr}$$
   where $P_j$ is the transition  operator associated to
$p(j,x,y)$. Hence
$$\E[U_n- U_{n-1}-G(x)\mid
\F_{n-1}]=\sum_{i=0}^{n-1}
[P_1G(X_{n-1}-X_i-x)-G(X_{n-1}-X_i-x)].$$
   Comparing with (4.1) and using
$$P_1G(z)-G(z)=-p(1,0,z), \qq z\in \Z^2,\eqno (4.3)$$ we see that
$$\E[ M_n- M_{n-1}\mid \F_{n-1}]=\E[ U_n- U_{n-1}+ \bar  B_n- \bar
B_{n-1}\mid \F_{n-1}]=0$$ as required.
\qed

The key to proving Theorem 1.1 in the $k=2$ case is the following
proposition.

\proclaim Proposition 4.2. We have
$$\E|\wt\beta_2(n,x)-\wt\beta_2(n,x')|^p
\leq c(p) (\log n)^{p}n |x-x'|^{p/3}
\eqno (4.4) $$ for each integer $p> 1$ and
$x,x'\in \Z^2/\sqrt n$ with $|x|,|x'|\leq 1$.

   Let
$$ W_2(n)=|B_2(n,x) |+|B_2(n,x') |,$$
$$ Y_2(n)=\max_{i\leq n}
\{| B_2(i,x)- B_2(i-1,x)|+| B_2(i,x')-B_2(i-1,x')|\},$$ and
$$\wt Z_2(n)=\sup_{j\leq n}|\wt U_2(j,x)-\wt U_2(j,x')|.$$

In the proof of Proposition 4.2 we will need the following three lemmas,
whose proofs are deferred until the next section.

\proclaim Lemma 4.3. For  any $x,x'$
  with $|x|,|x|'\leq \sqrt n$
$$\E  W_2(n)^p \leq c (\log n)^{p} n^p. \eqno (4.5)$$

\proclaim Lemma 4.4. For  any $x,x'$
  with $|x|,|x|'\leq \sqrt n$
$$\E  Y_2(n)^p\leq c n (\log n)^{p}. \eqno (4.6)$$

\proclaim Lemma 4.5. For any $x,x'$
  with $|x|,|x|'\leq \sqrt n$
$$ \E \wt Z_2(n)^p \leq c n^p \Bigl|\frac{x-x'}{\sqrt{n}}\Bigl|^{2 p/3},
\eqno (4.7)$$

\longproof{of Proposition 4.2} Converting from $\wt \beta$'s to $\wt B$'s,
estimate (4.4) for
$k=2$ is equivalent to
$$\E |\wt B_2(n,x)-\wt B_2(n,x')|^p\leq c(p) (\log n)^{p}  n^{p+1}
\Big(\frac{|x-x'|}{\sqrt n}\Big)^{p/3} \eqno (4.8)$$ for $x,x' \in \Z^{ 2}$
with
$|x|,|x'|\leq \sqrt n$. We want to apply Proposition 3.2.  We fix an $n$.  We
use the notation $f^{ +}( x)=\max(f( x),0 )$, $f^{ -}( x)=\max(-f( x),0 )$ so
that $f( x)=f^{ +}( x)-f^{-}( x)$. Take for $i\leq n$
$$Q^{ 1}_{ i}= B_2(i,x) +(G^{ +}(x')+G^{ -}(x) )i, \qq Q^{ 2}_{ i}=
B_2(i,x')+(G^{ +}(x)+G^{ -}(x') )i,\eqno(4.9)$$ so that $Q^1$ and
$Q^2$ are increasing and
$Q_i=Q^1_i-Q^2_i=\wt B_2(i,x)-\wt B_2(i,x')$.  For $i\geq n$ and $j=1,2$,
set $Q^j_i=Q^j_n$. We set  $H_{ i}=\wt U_2(i,x)-\wt U_2(i,x')$. By
Proposition 4.1,
   $Q_{ i}+H_{ i}$ is a martingale.

{}From Proposition 2.1, Lemmas 4.3 and 4.4, and the fact that $|x|, |x'|\leq
\sqrt n$, we see that
$$\E[(Q_n^1+Q_n^2)^p]\leq c(\log n)^p n^p \eqno (4.10)$$ and
$$\E[(\max_{i\leq n}\{ [Q^1_i-Q^1_{i-1}]+[Q^2_i-Q^2_{i-1}]\})^p]
\leq cn (\log n)^p. \eqno (4.11)$$ Combining (4.10), (4.11), Lemma 4.5,
and the fact that
$\frac{1}{\sqrt n}\leq \frac{|x-x'|}{\sqrt n}\leq 2$ unless $x=x'$ with
Proposition 3.2, we obtain
$$\E \sup_{j\leq n}|\wt B_2(j,x)-\wt B_2(j,x')|^p\leq c(p) (\log n)^{p} n^{
p+1}
\Big(\frac{|x-x'|}{\sqrt n}\Big)^{p/3} \eqno (4.12)$$ for $x,x' \in \Z^{ 2}$
with
$|x|,|x'|\leq \sqrt n$, which implies (4.8).
  This is the bound we need.
\qed

\longproof{of Theorem 1.1, the $k=2$ case}
  Let $f:\R^2\to [0,\infty)$ be a nonnegative $C^\infty$  function with
support  in $\{y: \frac12\leq |y|\leq 1\}$, and with integral 1. Let
$f_\tau(x)=\tau^{-2} f(x/\tau)$. The gradient of $f_\tau$ is bounded by a
constant times
$\tau^{-3}$. Set
$\tau_n=n^{-\zeta/4}$. Then recalling (1.7),
$$\int_0^1
\int_0^t|f_{\tau_n}(X^n_t-X^n_s)-f_{\tau_n}(W^n_t-W^n_s)|  ds \, dt
\leq c
\tau_n^{-3} n^{-\zeta}\leq c n^{-\zeta/4}. \eqno (4.13)$$  We also have by
Lemma 2.3 that
$$\Biggl|\sum_{ x\in \Z^{ 2}/\sqrt n} f_{\tau_n}(x) G_{ n}(x\sqrt n)
{1\over n} -\sum_{ x\in \Z^{ 2}/\sqrt n} f_{\tau_n}(x)
\frac{1}{\pi}\log ( 1/|x|) {1\over n}\Biggr|\leq c n^{-\bar \delta}\eqno
(4.14)$$  and it is easy to see from the support properties of
$f_{\tau_n}(x)$ that
$$\Biggl|\int f_{\tau_n}(x) \tfrac{1}{\pi}
\log(1/|x|)\, dx-\sum_{ x\in \Z^{ 2}/\sqrt n} f_{\tau_n}(x)
\frac{1}{\pi}\log ( 1/|x|) {1\over n}
\Biggr|\leq c n^{-\bar \delta}.\eqno (4.15)$$ On the other hand,
recalling the notation $X^n_t= X_{[nt]}/\sqrt n$
$$\eqalignno{ \sum_{ x\in \Z^{ 2}/\sqrt n}
f_{\tau_n}(x)  B_2(n,\sqrt n x)
{1\over n^{ 2}} &= {1\over n^{ 2}}\sum_{ x\in \Z^{ 2}/\sqrt n}
f_{\tau_n}(x)
\sum_{0\leq i_{ 1}<i_{ 2}\leq n}
   1_{(X_{ i_{ 2}}=X_{ i_{ 1}}+ \sqrt n x)}  \cr& ={1\over n^{ 2}}
\sum_{0\leq i_{ 1}<i_{ 2}\leq n} f_{\tau_n}\Big({(X_{ i_{ 2}}-X_{ i_{ 1}}
\over \sqrt n}\Big)
    \cr& ={1\over n^{ 2}}\int_0^n \int_0^t f_{\tau_n}\Big({(X_{[t]}-X_{[s]}
\over \sqrt n}\Big) ds\, dt
   \cr & =\int_0^1 \int_0^t f_{\tau_n}(X^n_t-X^n_s) ds\, dt &(4.16)\cr}$$
so that
$$\sum_{ x\in \Z^{ 2}/\sqrt n} f_{\tau_n}(x) \wt
\beta_2(n,x)  {1\over n}=\int_0^1 \int_0^t f_{\tau_n}(X^n_t-X^n_s) ds\, dt
   -\sum_{ x\in \Z^{ 2}/\sqrt n} f_{\tau_n}(x) G_{ n}(x\sqrt n) {1\over n}
.$$   By \cite{\BKb},
$$\int f_{\tau_n}(x) \wt \gamma_2(1,x,n) \, dx =\int_0^1 \int_0^t
f_{\tau_n}(W^n_t-W^n_s) \, ds\, dt -\int f_{\tau_n}(x) \tfrac{1}{\pi}
\log(1/|x|)\, dx. \eqno (4.17)$$ (This conforms with the definition given in
Section 1 above; the definition in \cite{\BKb} is very slightly different and
would yield
$\int_0^1 \int_0^1 f_{\tau_n} (W^n_t-W^n_s) \, ds\, dt$ instead.)
Combining the above,
$$\Biggl|\sum_{ x\in \Z^{ 2}/\sqrt n} f_{\tau_n} (x) \wt \beta_2(n,x) \,
{1\over n}-
\int f_{\tau_n}(x) \wt \gamma_2(1,x,n)\, dx\Biggr|
=O(n^{-\zeta/4})+O(n^{-\ol \delta}). \eqno (4.18)$$
  Recall $\int f_{\tau_n}(x) dx=1$. Without loss of generality we may assume
$\zeta $ is small enough so that $\psi_{ n}=:\sum_{ x\in \Z^{ 2}/\sqrt n}
f_{\tau_n}(x) {1\over n}=1+O(n^{-\ol \delta} )$ for some $\ol \delta>0$. (If
$\zeta$ were too large, then $\tau_n$ would tend to 0 too quickly, and then
the above estimate for $\psi_n$ might not be valid. In general one only has
$\psi_n=1+O(n^{-1/2}\tau_n^{-3})$. )
  Jensen's inequality and estimates (4.4), (4.5) imply that
$$\eqalignno{\E \Biggl| \sum_{ x\in \Z^{ 2}/\sqrt n} f_{\tau_n}(x) \wt
\beta_2(n,x) {1\over n}- \psi_n \wt
\beta_2(n,0)\Biggr|^p &\leq \sum_{ x\in \Z^{ 2}/\sqrt n} f_{\tau_n}(x)
\E|\wt
\beta_2(n,x)-\wt
\beta_2(n,0)|^p  {1\over n}\cr &\leq c(p)(\log n)^{p}n (\tau_n)^{p/3}.&(
4.19)\cr}$$ If we take  $p$ big enough, then
$$\P\Big(\Bigl|\sum_{ x\in \Z^{ 2}/\sqrt n} f_{\tau_n}(x) \wt
\beta_2(n,x) {1\over n} -\psi_n \wt
\beta_2(n,0)\Bigr|
\geq n^{-\zeta /24}\Bigr)
\leq c\frac{(\log n)^{p}n (\tau_n)^{p/3}}{n^{-\zeta p/24}}
\leq \frac{c}{n^2}.$$ By Borel-Cantelli, we conclude that
$$\sum_{ x\in \Z^{ 2}/\sqrt n} f_{\tau_n}(x) \wt
\beta_2(n,x) {1\over n} -\psi_n\wt
\beta_2(n,0)=O(n^{-\zeta/24}),
\qq a.s. \eqno (4.20)$$ Using (4.10),
$$\E[|\wt B_2(n,0)^p|]\leq c(\log n)^p n^p,$$ so
$$\P(|\wt \beta_2(n,0)|>n^{\ol \delta/2})\leq \frac{c(\log n)^p}{n^{\ol
\delta p/2}},$$ and if we take $p$ large enough, Borel-Cantelli tells us that
$$\wt \beta_2(n,0)=O(n^{\ol \delta/2}), \qq a.s.$$ So then
$$\wt \beta_2(n,0)-\psi_n \wt \beta_2(n,0)=O(n^{-\ol \delta/2}), \qq
a.s.$$

A very similar argument to the above also shows that we have
$$\int f_{\tau_n}(x) \wt \gamma_2(1,x,n)\, dx-\wt
\gamma_2(1,0,n)=O(n^{-\zeta/24}),
\qq a.s. \eqno (4.21)$$ the analogue to estimate (4.4) is in
\cite{\BKb}.
   Combining, we conclude that
$$\wt \beta_2(n,0)-\wt
\gamma_2(1,0,n)=O(n^{-\zeta /24})+O(n^{-\ol \delta/2}),
\qq a.s.$$
\qed

\ms

\ni {\bf Remark 4.6.} To see the importance of renormalization, note that if
we also had the estimate (4.8) for
$ B_2(n,x)-B_2(n,x')$, this would imply that uniformly in $n$
$$ |G( x)-G( x')|\leq c( p) ( \log n)\, n^{1/p}
\Big(\frac{|x-x'|}{\sqrt n}\Big)^{1/3}$$ which is impossible if
$p>6$ and $n$ is sufficiently large.
\ms

\subsec{5. Proofs of Lemmas 4.3-4.5}

\longproof{of Lemma 4.3} We have
$$\eqalign{ \E\{&(B_2(n,x))^{ p}\}\cr &\leq
\E\{(\sum_{i,j=0}^n 1_{(X_j-X_i=x)})^{ p}\}
\cr &=\sum_{s\in \cal{S} }\,\sum_{0\leq i_{ 1}\leq \ldots i_{2p}\leq
n}\,\sum_{z_{ 1},\ldots,z_{ p}\in\Z^{ 2}}
\prod_{ j=1}^{ 2p} p(i_{ j} -i_{ j-1},z_{ s( j-1)} +x_{c( j-1) }, z_{ 
s( j)} +x_{c( j)
})
\cr}\eqno(5.1)$$ where $s$ runs over the set of maps $\cal{S}$ from
$\{ 1,\ldots,2p\}$ to $\{ 1,\ldots,p\}$ such that $s^{ -1}(j)=2$ for each
$1\leq j\leq p$,
   $c( j)=\sum_{ i=1}^{ j}1_{ \{ s( i)=s(j)\}}$ and $x_{1}=0, x_{2}=x$. Here
we use the conventions $i_{ 0}=0, z_{ 0}=0, c( 0)=0.$ Setting
$$  g_{ n}( x)=\sum_{i=0}^n p( i,0,x)\leq c\log n\eqno(5.2)$$
for $x\leq \sqrt n$ by Proposition 2.1, we can
bound (5.1) by
$$\eqalign{ \sum_{s\in \cal{S}}\, &
\sum_{z_{ 1},\ldots,z_{ p}\in\Z^{ 2}}\,
\prod_{ j=1}^{ 2p} g_{ n}(z_{ s( j)}+x_{c( j) }-z_{ s( j-1)}-x_{c( j-1) })\cr
&\leq c(\log n)^{ p}
\sum_{s\in \cal{S}}\sum_{z_{ 1},\ldots,z_{ p}\in\Z^{ 2}}
\prod_{ j\,:\,c( j)=1} g_{ n}(z_{ s( j)}+x_{c( j) }-z_{ s( 
j-1)}-x_{c( j-1) })\cr
&\leq c(n\log n)^{ p}.}\eqno(5.3)$$ Here we used the obvious fact that
$\sum_{x\in\Z^{ 2}}g_{ n}( x)=\sum_{x\in\Z^{ 2}}\sum_{i=0}^n p(
i,0,x)=n+1$.
\qed

\longproof{of Lemma 4.4} Let
$$C(m,x)=[B_2(m,x)-B_2(m-1,x)]=\sum_{i=0}^{m-1} 1_{(X_m=X_i+x)}.$$
We show
$$\E C(i,x)^p\leq c(\log n)^{p};\eqno(5.4)$$ using
$$Y_2(n)^p \leq \sum_{i=1}^n c(p)C(i,x)^p,$$
   we are then done.

But
$$\eqalign{& \E\{(\sum_{i=0}^{m-1} 1_{(X_m=X_i+x)})^{ p}\}\cr
&=p!\sum_{0\leq i_{ 1}\leq \ldots i_{p}\leq m-1}\sum_{z_{
1},\ldots,z_{ p},y\in\Z^{ 2}}
\prod_{ j=1}^{ p} p(i_{ j} -i_{ j-1},z_{ j-1}, z_{  j})1_{(y=z_{
j}+x)}p(m -i_{p},z_{p}, y)
\cr &=p!\sum_{0\leq i_{ 1}\leq \ldots i_{p}\leq m-1}\sum_{y\in\Z^{
2}}p(i_{ 1},0, y-x)
\prod_{ j=2}^{ p} p(i_{ j} -i_{ j-1},0, 0)p(m -i_{p},0, x)
\cr  &=p!\sum_{0\leq i_{ 1}\leq \ldots i_{p}\leq m-1}
\prod_{ j=2}^{ p} p(i_{ j} -i_{ j-1},0, 0)p(m -i_{p},0, x)
\cr&\leq c(\log m)^{ p}
\cr}$$ which is (5.4).
\qed

\longproof{of Lemma 4.5} We begin by  estimating
$$\E (1+|X_i|^2)^{-b/2}$$   with $0<b< 2$.

First
$$\E[(1+|X_i|^2)^{-b/2}; |X_i|=0]=1\cdot \P(X_i=0)\leq c/i\leq ci^{-b/2}.$$
Next, using (2.1)
$$\eqalign{\E[(1+|X_i|^2)^{-b/2}; 0<|X_i|<\sqrt i]&
\leq \sum_{\{x\in \Z^2, 0<|x|\leq \sqrt i\}} |x|^{-b} p(i,0,x)\cr &\leq
\frac{c}{i} \sum_{\{x\in \Z^2, 0<|x|\leq \sqrt i\}} |x|^{-b}\cr &=\frac{c}{i}
i^{1-b/2}=ci^{-b/2}.\cr}$$ Finally, 
$$\E[(1+|X_i|^2)^{-b/2}; \sqrt i\leq |X_i|\,]
\leq
  (1+ i)^{-b/2} \P(\sqrt i\leq |X_i|)
\leq ci^{-b/2}.\eqno (5.5)$$ We conclude that for any $0<b< 2$
$$\E[(1+|X_i|^2)^{-b/2}]\leq ci^{-b/2}.\eqno (5.6)$$

Using the estimate
$$|G(X_i+x)-G(X_i+y)|\leq \frac{c|x-y|^{ 2/3}}{(1+|X_i+x|^2)^{1/3}}
+\frac{c|x-y|^{ 2/3}}{(1+|X_i+y|^2)^{1/3}}$$ of Proposition 2.2, the fact
that  symmetry tells us that
$\E[(1+|X_i+x|^2)^{-1/3}]$ is largest when $x=0$, and the  estimate (5.6)
above, we obtain
$$\E\sum_{i=1}^n |G(X_i+x)-G(X_i+y)|\leq c|x-y|^{ 2/3}
\sum_{i=1}^n i^{-1/3}
\leq cn(|x-y|/\sqrt n)^{ 2/3}.$$  So by independence, using
$\bar X_i, \bar \E$ to denote an independent copy of $X_i$ and its
expectation operator,
$$\eqalignno{&\E\Big[ \sum_{i=m+1}^n |G(X_i+x)-G(X_i+y)|\ \mid
\F_m\Big]&(5.7)\cr &
\leq \bar \E \sum_{i=1}^n|G(\bar X_i+X_m+x)-G(\bar X_i+X_m+y)|\leq
cn\Bigl(\frac{|x-y|}{\sqrt n}\Big) ^{2/3}.\cr}$$  If $|x|\leq \sqrt n$, by
Proposition 2.1 and Doob's inequality
$$\eqalignno{\E \sup_{i\leq n} &|G(X_i+x)|^p
\leq c\E \sup_{i\leq n}(1+\log^+|X_i+x|)^p&(5.8)\cr &\leq c(\log n)^p
\P(\sup_{i\leq n} |X_i+x|\leq n) +c\sum_{m=[\log n]}^\infty m^p
\P(e^m\leq \sup_{i\leq n} |X_i+x|\leq e^{m+1})\cr &\leq c(\log n)^p
+c\sum_{m=[\log n]}^\infty  m^p \frac{\E |X_n+x|^2}{e^{2m}}\cr &\leq
c(\log n)^p.\cr}$$ If $x\ne y$, then
$|x-y|\geq 1$ and (5.8) then implies that
$$\E \sup_{i\leq n} |G(X_i+x)-G(X_i+y)|^p\leq c(\log n)^p
\leq c\Big(n\Big(\frac{|x-y|}{\sqrt n}\Big)^{2/3}\Big)^p.$$
  Using Proposition 3.1 with $A_j=\sum_{i=1}^{j\land n}
|G(X_i+x)-G(X_i+y)|$, if $|x|, |y|\leq
\sqrt n$, then
$$\Bigl\|\sum_{i=1}^n |G(X_i+x)-G(X_i+y)|\,\Bigr\|_{ p}
\leq cn(|x-y|/\sqrt n)^{ 2/3}.\eqno(5.9)$$  Replacing $x$ and $y$ by
$-x$ and $-x'$, resp., and using the fact that $\sum_{i=1}^n G(X_i-x)$ is
equal in law to
$\sum_{i=0}^{n-1} G(X_n-X_i-x)$ yields the
$L^p$ estimate that we want.
\qed

\subsec{6. The martingale connection:  $k>2$}

Let $\wt B_{ 1,m}(j,x)=j$ and for $k\geq 2$ define
$$\wt B_{ k,m}(j,x)=\sum_{A\subset\{2, \ldots,k\}} (-1)^{|A|}
\Big(\prod_{i\in A} G_{ m}(x_i)\Big) B_{k-|A|}(j,x_{A^c}).
\eqno (6.1)$$ Note that $\wt B_{ k,n}(n,x)=\wt B_{ k}(n,x)$.

If $x=(x_2, \ldots, x_{k-1}, x_{k})$, let $x_{k^c}=(x_2, \ldots, x_{k-1})$.
\proclaim Proposition 6.1. Let $k>2$. If
$$\wt  U_{ k,m}(n,x)=\sum_{i=1}^{n} G_m(X_n-X_i-x_{ k}) [\wt B_{
k-1,m}(i,x_{ k^{ c}})-\wt B_{ k-1,m}(i-1,x_{ k^{ c}})],$$ then for each
$m$
$$M_{ n,m}=\wt  U_{ k,m}(n,x)+
\wt B_{ k,m}(n,x),\hskip.4in n=0,1,2,\ldots$$ is a martingale with
$M_{ 0,m}=0$.

\proof We will show that for each $k$
$$ N_{ k,m}( n)=U_{ k,m}(n,x)+B_{ k,m}(n,x)-G_{ m}(x_k)B_{
k-1,m}(n,x_{ k^{ c}}),\hskip.4in n=0,1,2,\ldots$$ is a martingale where
$$\eqalign{U_{ k,m}(n,x)&=\sum_{i=1}^{n} [G(X_n-X_i-x_{ k})-G(e_{
1}\sqrt{m})][ B_{ k-1}(i,x_{ k^{ c}})-B_{ k-1}(i-1,x_{ k^{ c}})]\cr &=
\sum_{i=1}^{n} G(X_n-X_i-x_{ k})[ B_{ k-1}(i,x_{ k^{ c}})-B_{ k-1}(i-1,x_{
k^{ c}})]\cr &\qquad\qquad\qquad\qquad-G(e_{1}\sqrt{m})B_{ k-1}(n,x_{
k^{ c}}).\cr}$$

This will prove the proposition since, with the notation $D_{ k}=D\cup\{
k\}$,
$$\eqalign{&\wt  U_{ k,m}(n,x)+
\wt B_{ k,m}(n,x)\cr &=
\sum_{i=1}^{n} G_m(X_n-X_i-x_{ k})
\cr&  \sum_{D\subset\{2, \ldots,k-1\}}(-1)^{|D|}
\Big(\prod_{l\in D} G_{ m}(x_l)\Big) [B_{ k-|D_{ k}|,m}(i,x_{ D_{ k}^{ c}})-
B_{ k-|D_{ k}|,m}(i-1,x_{ D_{ k}^{ c}})]\cr & +\sum_{A\subset\{2, \ldots,k\}}
(-1)^{|A|}
\Big(\prod_{i\in A} G_{ m}(x_i)\Big) B_{k-|A|}(n,x_{A^c})\cr &
=\sum_{D\subset\{2, \ldots,k-1\}}(-1)^{|D|}
\Big(\prod_{l\in D} G_{ m}(x_l)\Big)\cr & \hskip.5in [U_{ k-|D^{
c}|,m}(n,x_{D^{ c} })+B_{ k-|D^{ c}|,m}(n,x_{D^{ c} })-G_{ m}(x_k)B_{ k-|D^{
c}|-1,m}(n,x_{D_{ k}^{ c} })].\cr}$$

If we set
$$\bar U_k(n,x)=\sum_{i=1}^{n} G(X_n-X_i-x_{ k})[ B_{ k-1}(i,x_{ k^{
c}})-B_{ k-1}(i-1,x_{ k^{ c}})]$$ we have that
$$ N_{ k,m}( n)=\bar U_k(n,x)+B_k(n,x)-G(x_k)B_{ k-1}(n,x_{ k^{ c}}).$$
Abbreviating $\bar U_n=\bar U_k(n,x)$ and
$\bar B_n=B_k(n,x)-G(x_k)B_{ k-1}(n,x_{ k^{ c}})$, we have that
$$ N_{ k,m}( n)=\bar U_n+\bar B_n.$$ Setting
$$H_{ i}=B_{ k-1}(i,x_{ k^{ c}})-B_{ k-1}(i-1,x_{ k^{ c}}) =\sum_{0\leq i_{
1}<i_{ 2}<\ldots<i_{ k-1}= i}\,
\prod_{j=2}^{k-1} 1_{(X_{ i_{ j}}=X_{ i_{ j-1}}+x_j)}$$ we have
$$\bar B_n-\bar B_{n-1}=\sum_{i=1}^{n-1} 1_{(X_n=X_i+x_k)}H_{
i}-G(x_k)H_{ n}.$$  So using (4.1)
$$\eqalign{\E[\bar B_n &-\bar B_{n-1}+G(x_k)H_{ n}\mid
\F_{n-1}]\cr  &=\sum_{i=1}^{n-1} p(1,0,X_{n-1}-X_i-x_k)H_{
i}.\cr}\eqno(6.2)$$
    From the definition of $\bar U_n$ we have
$$\eqalign{\bar U_n-\bar U_{n-1}&=G(x_k)H_{ n}\cr &\qq
+\sum_{i=1}^{n-1} [G(X_n-X_i-x_k)-G(X_{n-1}-X_i-x_k)]H_{ i}.\cr}$$
Recalling (4.2)
$$\eqalignno{\E[\bar U_n-\bar U_{n-1}&-G_{ n}(x_k)H_{ n}\mid
\F_{n-1}]&\cr &\qq =\sum_{i=1}^{n-1}
[P_1G(X_{n-1}-X_i-x_k)-G(X_{n-1}-X_i-x_k)]H_{ i}.\cr}$$ Comparing with
(6.2) and using (4.3)
$$P_1G(x)-G(x)=-p(1,0,x),$$ we see that
$$\E[\bar U_n-\bar U_{n-1}+\bar B_n-\bar B_{n-1}\mid
\F_{n-1}]=0$$ as required.
\qed

\ms
\ni{\bf Remark.} The statement of Proposition 6.1 is not an exact analogue
of that of Proposition 4.1. Consider the summands in the definition of
$\wt U_{ k,m}(n,x)$:
$$G_m(X_n-X_i-x_k) [\wt B_{k-1,m}(i,x_{k^c})-\wt
B_{k-1,m}(i-1,x_{k^c})].\eqno (6.3)$$ When $k=2$ and $i=n$, this is
nonrandom, whereas this is not the case when $k>2$ and $i=n$.  On the
other hand, recalling that $B_{k-1}(i,x)=0$ if
$i>k-1$, it is natural to define $B_{k-1}(-1,x)$ to be 0. It is also natural to
define $\wt B_{ 1,m}(i,x)=i$ for $i\geq 0$. Then (6.3) will be 0 if
$i=0$ for all $k\geq 2$, but the $i=0$ term in the statement of Proposition
4.1 is not zero.
\ms

\subsec{7. The case of general $k$}

The key to proving Theorem 1.1 for the case of general $k$ is the following
proposition.

\proclaim Proposition 7.1. For any $k\geq 2$ we have
$$\E|\wt\beta_k(n,x)-\wt\beta_k(n,x')|^p
\leq c(p) (\log n)^{p( k-1)}n^k |x-x'|^{8^{ -k}p}
\eqno (7.1) $$ for each integer $p> 1$ and $x,x'\in (\Z^2)^{ k-1}/\sqrt n$
with $|x|,|x'|\leq 1$.

   Let
$$W_k(n)=| B_k(n,x)|+| B_k(n,x')|,$$
$$\wt W_k(n)=|\wt B_k(n,x)|+|\wt B_k(n,x')|,$$
$$ Y_k(n)=\max_{i\leq n}
\{| B_k(i,x)- B_k(i-1,x)|+| B_k(i,x')- B_k(i-1,x')|\},$$
$$\wt Y_k(n)=\max_{i\leq n}
\{|\wt B_{ k,n}(i,x)-\wt B_{ k,n}(i-1,x)|+|\wt B_{ k,n}(i,x')-\wt B_{
k,n}(i-1,x')|\},$$
  and
$$\wt Z_k(n)=\sup_{j\leq n}|\wt U_{ k,n}(j,x)-\wt U_{ k,n}(j,x')|.$$

In the proof of Proposition 7.1 we will need the following three lemmas,
whose proofs are deferred until the next two sections.

\proclaim Lemma 7.2. For  any $x,x'$
  with $|x|,|x|'\leq \sqrt{n}$
$$\E  W_k(n)^p \leq c (\log n)^{p( k-1)} n^p. \eqno (7.2)$$ and
$$\E \wt W_k(n)^p \leq c (\log n)^{p( k-1)} n^p. \eqno (7.3)$$

\proclaim Lemma 7.3. For  any $x,x'$
  with $|x|,|x|'\leq \sqrt{n}$
$$\E Y_k(n)^p\leq c n (\log n)^{p( k-1)}. \eqno (7.4)$$  and
$$\E \wt Y_k(n)^p\leq c n (\log n)^{p( k-1)}. \eqno (7.5)$$

\proclaim Lemma 7.4. For any $x,x'$
  with $|x|,|x|'\leq \sqrt{n}$
$$ \E \wt Z_k(n)^p \leq c (\log n)^{p( k-1)} n^{ p+k}
\Big(\frac{|x-x'|}{\sqrt n}\Big)^{8^{ -k}2p} ,
\eqno (7.6)$$

\longproof{of Proposition 7.1} Converting from $\wt \beta$'s to $\wt B$'s,
estimate (7.1) is equivalent to
$$\E |\wt B_{ k,n}(n,x)-\wt B_{ k,n}(n,x')|^p\leq c(p) (\log n)^{p( k-1)}
n^{ p+k}
\Big(\frac{|x-x'|}{\sqrt n}\Big)^{8^{ -k}p} \eqno (7.7)$$ for $x,x' \in
(\Z^2)^{ k-1}$ with
$|x|,|x'|\leq \sqrt{n}$. We want to apply Proposition 3.2.  We fix an $n$.  Let
$${\cal A}_{ +}( x)=\bigg\{A\subset\{2, \ldots,k\}\,\big|
\, ( -1)^{ |A|}\prod_{i\in A} G_{ n}(x_i)> 0\bigg\}$$
$$\wt B_{ k,n,+}(j,x)=\sum_{A\in {\cal A}_{ +}( x)}( -1)^{ |A|}
\Big(\prod_{i\in A} G_{ n}(x_i)\Big) B_{k-|A|}(j,x_{A^c}).
\eqno (7.8)
$$
$$\wt B_{ k,n,-}(j,x)=\sum_{A\in {\cal A}^{ c}_{ +}( x)}( -1)^{ |A|}
\Big(\prod_{i\in A} G_{ n}(x_i)\Big) B_{k-|A|}(j,x_{A^c}).
\eqno (7.9)
$$  For $i\leq n$ set
$$Q^{ 1}_{ i}= \wt B_{ k,n,+}(i,x)-\wt B_{ k,n,-}(i,x'), \qq Q^{ 2}_{ i}=
\wt B_{ k,n,+}(i,x')-\wt B_{ k,n,-}(i,x),$$ so that $Q^1$ and $Q^2$ are
increasing and
$Q_i=Q^1_i-Q^2_i=\wt B_{ k,n}(i,x)-\wt B_{ k,n}(i,x')$.  For
$i\geq n$ set $Q^j_i=Q^j_n$, $j=1,2$. We set  $H_{ i}=\wt U_{ k,n}(i,x)-\wt
U_{ k,n}(i,x')$. By Proposition 6.1,
   $Q_{ i}+H_{ i}$ is a martingale. Using  Lemmas 7.2-7.4 and Proposition 2.1
to bound the right hand side of  (3.5) in Proposition 3.2 and using the fact
that
$\frac{1}{\sqrt n}\leq \frac{|x-x'|}{\sqrt n}\leq 2$ unless $x=x'$ we obtain
$$\E \sup_{j\leq n}|\wt B_{ k,n}(j,x)-\wt B_{ k,n}(j,x')|^p
\leq c(p) (\log n)^{p( k-1)} n^{p+k}
\Big(\frac{|x-x'|}{\sqrt n}\Big)^{8^{ -k}p} \eqno (7.10)$$  for $x,x' \in
(\Z^2)^{k-1}$
  with
$|x|,|x'|\leq \sqrt{n}$, which implies (7.7).
  This is  the bound we need.
\qed

\longproof{of Theorem 1.1, the general case}  The proof is quite similar to
the $k=2$ case.
  Let $f:\R^2\to [0,\infty)$ be a nonnegative $C^\infty$ function with
support  in $\{y: \frac12\leq |y|\leq 1\}$, and with integral 1. Let
$f_\tau(x)=\tau^{-2} f(x/\tau)$.
  Set
$\tau_n=n^{-\zeta/4k}$.  Set
$g_{ n}(f)=\sum_{ x\in \Z^{ 2}/\sqrt n} f(x) G_{ n}(x\sqrt n)  {1 \over n}$
and
$ l(f)=\frac{1}{\pi}\int f(x) \log (   1/|x|) dx$. As in (4.14), (4.15)
$$  |g_{ n}(f_{\tau_n})-
\frac{1}{\pi}l(f_{\tau_n})| \leq c (\tau_n\sqrt n)^{-\rho}.
\eqno(7.11)$$  Using  (1.1) and setting $F_{\tau_n}(x_{2},\ldots,x_{ k}
)=\prod_{ i=2}^{ k} f_{\tau_n}( x_{ i})$ we have
$$\eqalign{& \int   F_{\tau_n}( x)
\wt\gamma_{k}(1,x)\,dx=\sum_{j=0}^{ k-1}\choose{k-1}{j}   ( -1)^{j}
\Big( l(f_{\tau_n})\Big)^{j}
\cr & \qquad\qquad\qquad\qquad\qquad\times\int_{ 0\leq t_{
1}\leq\cdots\leq t_{ k-j}<1}
\prod_{ i=2}^{ k-j} f_{\tau_n}(W_{ t_{ i}}-W_{ t_{ i-1}})dt_{ 1}\cdots dt_{
k-j}.
\cr}\eqno(7.12)$$
  On the other hand, as in (4.16), it is easily checked that we have
$$\eqalign{& \sum_{ x\in \Z^{ 2( k-1)}/\sqrt n}  F_{\tau_n}( x)
\wt\beta_{k}(n,x)\,{1 \over n^{ k-1}}=\sum_{j=0}^{ k-1}\choose{k-1}{j}   (
-1)^{j}
\Big( g_{ n}(f_{\tau_n})\Big)^{j}
\cr & \qquad\qquad\qquad\qquad\qquad\times\int_{ 0\leq t_{
1}\leq\cdots\leq t_{ k-j}<1}
\prod_{ i=2}^{ k-j} f_{\tau_n}(X^{ n}_{ t_{ i}}-X^{ n}_{ t_{ i-1}})dt_{
1}\cdots dt_{ k-j}.
\cr}\eqno(7.13)$$

Since  the gradient of $f_\tau$ is bounded by a constant times
$\tau^{-3}$
$$\eqalign{& \int_{ 0\leq t_{ 1}\leq\cdots\leq t_{ k-j}<1}
\Big|\prod_{ i=2}^{ k-j} f_{\tau_n}(W^{ n}_{ t_{ i}}-W^{ n}_{ t_{
i-1}})-\prod_{ i=2}^{ k-j} f_{\tau_n}(X^{ n}_{ t_{ i}}-X^{ n}_{ t_{ i-1}})\Big|
\,dt_{ 1}\cdots dt_{ k-j}
\cr & \qquad\qquad\qquad\qquad\qquad\qquad\qquad\leq c
\tau_n^{-3-2(k-j-2 )} n^{-\zeta}\leq c n^{-\zeta/4}.
\cr}\eqno(7.14)$$

Combining  (7.11), (7.14), and the fact that both $| g_{ n}(f_{\tau_n}) |$ and
$| l(f_{\tau_n}) |$ are bounded by $c\log n$ we see that
$$\eqalign{& \Big|\int   F_{\tau_n}( x)
\wt\gamma_{k}(t,x,n)\,dx-\sum_{ x\in \Z^{ 2( k-1)}/\sqrt n}  F_{\tau_n}( x)
\wt\beta_{k}(t,x)\,{1 \over n^{ k-1}}\Big|
\cr & \qquad\qquad\qquad\qquad\qquad\leq c
\tau_n^{-2(k-2 )} (\tau_n\sqrt n)^{-\rho}+ c ( \log n)^{ k-1}n^{-\zeta/4}
\leq c n^{-\zeta/8}
\cr}\eqno (7.15)$$ if we take $\zeta>0$ sufficiently small.

  Since $\int F_{\tau_n}(x) dx=1$, we have $\psi_{k, n}=:\sum_{ x\in
\Z^{ 2( k-1)}/\sqrt n} F_{\tau_n}(x) {1\over n^{ k-1}}=1+O(n^{-\ol
\delta} )$, provided we assume, as we may without loss of generality, that
$\zeta$ is sufficiently small.
  Jensen's inequality and estimates (7.1), (7.2) imply that
$$\eqalignno{&\E \Bigl| \sum_{ x\in \Z^{ 2( k-1)}/\sqrt n} F_{\tau_n}(x)
\wt
\beta_k(n,x) {1 \over n^{ k-1}}- \psi_{k,n}\wt
\beta_k(n,0)\Bigr|^p \cr &\leq \sum_{ x\in \Z^{ 2( k-1)}/\sqrt n}
F_{\tau_n}(x)
\E|\wt
\beta_k(n,x)-\wt
\beta_k(n,0)|^p  {1 \over n^{ k-1}}\cr &\leq c(p)(\log n)^{p( k-1)}n^2
(\tau_n)^{8^{ -k}p}.&(7.16)\cr}$$  If we take  $p$ big enough, then
$$\eqalign{&\P\Big(\Bigl|\sum_{ x\in \Z^{ 2( k-1)}/\sqrt n} F_{\tau_n}(x)
\wt
\beta_k(n,x) {1 \over n^{ k-1}} -\psi_{k,n}\wt
\beta_k(n,0)\Bigr|
\geq n^{-8^{ -k}\zeta/8k}\Bigr)\cr &
\leq c\frac{(\log n)^{p( k-1)}n^2 (\tau_n)^{8^{ -k}p}}{n^{-8^{
-k}p\zeta/8k}}
\leq \frac{c}{n^2}.\cr}$$ By Borel-Cantelli, we conclude that
$$\sum_{ x\in \Z^{ 2( k-1)}/\sqrt n} F_{\tau_n}(x) \wt \beta_k(n,x) {1
\over n^{ k-1}} -\psi_{k,n}\wt
\beta_k(n,0)=O(n^{-8^{ -k}\zeta /8k}),
\qq a.s. \eqno (7.17)$$ We use (7.3) and the same argument as in the
$k=2$ case to show
$$\wt \beta_k(n,0)-\psi_{k,n}\wt \beta_k(n,0)=O(n^{-\ol \delta/2}), \qq
a.s.
$$ for some $\ol \delta>0$. This with (7.17) yields
$$\sum_{ x\in \Z^{ 2( k-1)}/\sqrt n} F_{\tau_n}(x) \wt \beta_k(n,x) {1
\over n^{ k-1}} -\wt
\beta_k(n,0)=O(n^{-8^{ -k}\zeta /8k})+O(n^{-\ol \delta/2}),
\qq a.s. \eqno (7.18)$$

A similar argument shows that we have (7.18) holding with the
   $\wt \beta_k(n,x)$ replaced by
$\wt \gamma_k(1,x,n)$; the analogue to estimate (7.1) is in
\cite{\BKb}.
   Combining, we conclude that
$$\wt \beta_k(n,0)-\wt
\gamma_k(1,0,n)=O(n^{-8^{ -k}\zeta /8k})+O(n^{-\ol \delta/2}),
\qq a.s.$$
\qed

\subsec{8. Proofs of Lemmas 7.2-7.3}

These are again similar to the $k=2$ case.

\longproof{of Lemma 7.2} Using Proposition 2.1 it suffices to prove (7.2) for
all $k$.

We have
$$\eqalignno{ \E\{(B_k&(n,x))^{ m}\}&(8.1)\cr &\leq
\E\{(\sum_{i_{1},\ldots,i_{ k}=0}^n
\prod_{j=2}^{k} 1_{(X_{ i_{ j}}=X_{ i_{ j-1}}+x_j)})^{ m}\}
\cr &=\sum_{s\in S(k,m ) }\,\sum_{0\leq i_{ 1}\leq \ldots i_{km}\leq
n}\sum_{z_{ 1},\ldots,z_{ m}\in\Z^{ 2}}
\prod_{ j=1}^{ km} p(i_{ j} -i_{ j-1},z_{ s( j-1)} +\bar x_{c( j-1) 
}, z_{ s( j)}
+\bar  x_{c( j) })
\cr}$$ where $s$ runs over the set of maps $S(k,m ) $ from
$\{ 1,\ldots,km\}$ to $\{ 1,\ldots,m\}$ such that $s^{ -1}(j)=k$ for each
$1\leq j\leq m$,
   $c( j)=\sum_{ i=1}^{ j}1_{ \{ s( i)=s(j)\}}$,
$\bar  x_{j}=\sum_{ l=2}^{ j}x_{l}$
and
$\bar  x_{1}=0$. Here we use the conventions $i_{ 0}=0, z_{ 0}=0, c( 0)=0.$
Setting
$$  g_{ n}( x)=\sum_{i=0}^n p( i,0,x)\leq c\log n\eqno(8.2)$$ we can
bound (8.1) by
$$\eqalign{ \sum_{s\in S(k,m )}&\,
\sum_{z_{ 1},\ldots,z_{ m}\in\Z^{ 2}}\,
\prod_{ j=1}^{ km} g_{ n}(z_{ s( j)}+\bar  x_{c( j) }-z_{ s( j-1)}-\bar  x_{c(
j-1) })\cr &\leq c(\log n)^{( k-1) m}
\sum_{s\in \cal{S}}\,\sum_{z_{ 1},\ldots,z_{ k}\in\Z^{ 2}}\,
\prod_{ j\,:\,c( j)=1} g_{ n}(z_{ s( j)}+\bar  x_{c( j) }-z_{ s( 
j-1)}-\bar  x_{c(
j-1) })\cr &\leq cn^{ m}(\log n)^{( k-1) m}.}\eqno(8.3)$$ Here we used the
obvious fact that $\sum_{x\in\Z^{ 2}}g_{ n}( x)=\sum_{x\in\Z^{
2}}\sum_{i=0}^n p( i,0,x)=n+1$.
\qed

\longproof{of Lemma 7.3} Using Proposition 2.1 it suffices to prove (7.4) for
all $k$.

Let
$$C(n,x)=[B_k(n,x)-B_k(n-1,x)]=\sum_{0\leq i_{ 1}<i_{ 2}<\ldots<i_{ k}=n}
\prod_{j=2}^{k} 1_{(X_{ i_{ j}}=X_{ i_{ j-1}}+x_j)}.$$ If we show
$$\E C_k(i,x)^p\leq c(\log n)^{( k-1)p},\qquad i\leq n,\eqno(8.4)$$ then
using
$$Y_k(n)^p \leq \sum_{i=1}^n c(p)C_k(i,x)^p,$$
   we are done.

But
$$\eqalign{ \E\Big\{\Big(&\sum_{0\leq i_{ 1}<i_{ 2}<\ldots<i_{ k}=n}
\, \prod_{j=2}^{k} 1_{(X_{ i_{ j}}=X_{ i_{ j-1}}+x_j)}\Big)^{ m}\Big\}\cr
&=\sum_{s\in S(k-1,m ) }\,\sum_{(i_{ 1}, \ldots, i_{(k-1)m})\in D( s) }\,
\sum_{z_{ 1},\ldots,z_{ m},y\in\Z^{ 2}}\cr &
\qquad\prod_{ j=1}^{ (k-1)m} p(i_{ j} -i_{ j-1},z_{ s( j-1)} +\bar  x_{c( j-1)
}, z_{ s( j)} +\bar  x_{c( j) })
\cr &\qquad\qquad \Big( \prod_{ j\,:\,c( j)=k-1}1_{(y=z_{ s( j)}+\bar  x_{
k})}\Big) p(n -i_{(k-1)m}+\bar  x_{ k-1},z_{ s((k-1)m)},y)\cr}\eqno(8.5)$$
where
$S(k-1,m ), c( j), \bar  x_{j}$ are defined in the last section and for each
$s\in S(k-1,m )$
$$\eqalign{ D( s)=\{ (i_{ 1}, \ldots, i_{(k-1)m})\,:\,0&\leq i_{ 1}\leq
\ldots\leq i_{(k-1)m}<n
\,  \hbox{and}\, \cr &  i_{j-1}<i_j \hbox{ whenever } s( j-1)=s( j)\}.\cr}
$$ We can then see that (8.5) equals
$$
\eqalign{\sum_{s\in S(k-1,m ) }\,&\sum_{(i_{ 1}, \ldots, i_{(k-1)m})\in D(
s) }\sum_{y\in\Z^{ 2}}\cr &
\qquad p(i_{1},0,y-\bar  x_{ k})\prod_{ j=2}^{ (k-1)m} p(i_{ j}
-i_{ j-1},\bar  x_{c( j-1) },\bar   x_{c( j) })\,  p(n -i_{(k-1)m},0,x_{ k})\cr
&\leq c(\log n)^{(k-1)m}
\cr}$$ which is (7.4).
\qed

\subsec{9. Proof of Lemma 7.4}

This proof is substantially different from the proof of Lemma 4.5.

\def\didi{{[\Delta_{ i}\wt B_{ k-1,n}(i,x_{ k^{ c}})-\Delta_{ i}\wt B_{
k-1,n}(i,x'_{ k^{ c}})]}}
\def\dido{{[\Delta_{ i}\wt B_{ k-1,n}(i,x_{ k^{ c}})]}}

\longproof{of Lemma 7.4} We use induction on $k$. We already know (7.6)
for $k=2$. Thus assume (7.6) has been proved  with $k$ replaced by
$i$ for all $i\leq k-1$. Then as explained above in the proof of Proposition
7.1 we will have that (7.10) holds  with $k$ replaced by $i$ for all $i\leq
k-1$.

We will show that
$$\E[|\wt U_{k,n}(m,x)-\wt U_{k,n}(m,x)|^p]
\leq c(\log n)^{(k-1)p} n^{p+k-1}\Big(\frac{|x-x'|}{\sqrt n}\Big)^{8^{-k}2p}
\eqno (9.1)$$
for $m\leq n$. This and the inequality
$$\E[\max_{m\leq n} |C_m|^p]\leq \E \sum_{m=1}^n |C_m|^p
=\sum_{m=1}^n \E[|C_m|^p]\leq n\max_{m\leq n} \E[|C_m|^p]$$ yields
(7.6).

Abbreviating $\Delta_{ i}\wt B_{ k-1,n}(i,x_{ k^{ c}}):=\wt B_{ k-1,n}(i,x_{
k^{ c}})-\wt B_{ k-1,n}(i-1,x_{ k^{ c}})$, where $x_{k^c}$ is the same as
$(x_2, \ldots, x_{k-1})$, we have
$$\eqalignno{\|&\wt U_{ k,n}(m,x)-\wt U_{ k,n}(m,x')\|_{ p}&(9.2)\cr &
\leq \|\sum_{i=1}^{m} \Big(G(X_m-X_i-x_{ k})-G(X_m-X_i-x'_{
k})\Big)\,\Big(\Delta_{ i}\wt B_{ k-1,n}(i,x_{ k^{ c}})\Big)\|_{p}\cr &
\qq +\|\sum_{i=1}^{m} \Big(G(X_m-X_i-x'_{ k})-G(e_{
1}\sqrt{n})\Big)\,\Big(\Delta_{ i}\wt B_{ k-1,n}(i,x_{ k^{ c}}) -\Delta_{ i}\wt
B_{ k-1,n}(i,x'_{ k^{c}})\Big)\|_{p}.
\cr}$$
  Then with $m\leq n$
$$\eqalign{\|\sum_{i=1}^{m}& \Big(G(X_m-X_i-x_{ k})-G(X_m-X_i-x'_{
k})\Big)\,\Big(\Delta_{ i}\wt B_{ k-1,n}(i,x_{ k^{ c}})\Big)\|_{p}\cr &
\leq \|\sum_{i=1}^{m} |G(X_m-X_i-x_{ k})-G(X_m-X_i-x'_{
k})|\,\|_{2p}\,\|\wt Y_{ k-1}(m)\|_{2p}\cr & \leq c
n\Bigl|\frac{x-x'}{\sqrt{n}}\Bigr|^{2/3 }n^{ 1/2p} (\log n)^{( k-2)}.
\cr}\eqno (9.3)$$  by (5.9) and (7.5).

After interchanging  $x'$ and $x$ for convenience it remains to bound
$$\eqalign{&\Bigl\|\sum_{i=1}^{m} \Big(G(X_m-X_i-x_{
k})-G(e_{1}\sqrt{n})\Big)\,\Big(\Delta_{ i}\wt B_{ k-1,n}(i,x_{ k^{ c}})
-\Delta_{ i}\wt B_{ k-1,n}(i,x'_{ k^{c}})\Big)\Bigr\|_{p}\cr &
\leq \Bigl\|\sum_{i=1}^{m} G(X_m-X_i-x_{
k})\,\Big(\Delta_{ i}\wt B_{ k-1,n}(i,x_{ k^{ c}})
-\Delta_{ i}\wt B_{ k-1,n}(i,x'_{ k^{c}})\Big)\Bigr\|_{p}\cr &
+ \Bigl\|G(\sqrt n e_1)  [\wt B_{ k-1,n}(m,x_{ k^{ c}})-\wt B_{ k-1,n}(m,x'_{
k^{ c}})]\Bigr\|_{p}
\cr}\eqno (9.4)$$

   Using Proposition 2.1 and  our inductive hypothesis concerning
  (7.10) we see that 
$$\eqalign{&\Bigl\|G(\sqrt n e_1)  [\wt B_{ k-1,n}(m,x_{ k^{ c}})-\wt B_{
k-1,n}(m,x'_{ k^{ c}})]\Bigr\|_{p}\cr & \leq c (\log n)^{( k-1)} n^{1+(
k-1)/p}
\Big(\frac{|x-x'|}{\sqrt n}\Big)^{8^{ -k+1}}.\cr}\eqno (9.5)$$  To complete
the proof of (9.1) it therefore suffices to show that
$$\eqalign{\|\sum_{i=1}^{m}& G(X_m-X_i-x))\didi\|_{ p}\cr &\leq c (\log
n)^{( k-1)} n^{1+( k-1)/p}
\Big(\frac{|x-x'|}{\sqrt n}\Big)^{8^{ -k}2}.
\cr}\eqno (9.6)$$
\ms

By Propositions 2.1 and 2.3, $G(x)$ is bounded above (but $G(x)\to
-\infty$ as $|x|\to \infty$). Let
$$J(x)=G(x)\vee (-9\log m), \qq H(x)=G(x)-J(x).$$

Let $K_i=J(X_m-X_i-x)$ for $i=0, \ldots, m$ and let
$K_i=J(X_m-x)$ for $i<0$.
  Let $B$ be a small positive real  to be chosen later and let
$$L_i= \frac{K_i+\cdots +K_{i-Bm}}{Bm}.$$  We see that
$$\eqalignno{\Bigl| \sum_{i=1}^{m}& (K_i-L_i) \dido \Bigr| \cr &
\leq
\wt Y_{ k-1}( n)\sum_{i=1}^{n} |K_i-L_i|. &(9.7)\cr}$$

Since $J$ is bounded in absolute value by $c\log m$, the
  same is true for
$K_i$ and $L_i$ for any $i$, i.e.
$$\sup_{ i}|K_i-L_i| \leq c\log m. \eqno(9.8)$$
Note that $L_i$ and $K_i$ are independent of $\F_h$ for $i\geq h+Bm$,
and thus
$$\E\Bigl[ \sum_{i=h}^{m} |K_i-L_i| \mid \F_h\Big]
\leq \E \sum_{i=h+2Bm}^{m} |K_i-L_i|+ cBm\log m. \eqno (9.9)$$  Now by
Proposition 2.2
$$\eqalign{\E |J(X_m-X_i-x)&-J(X_m-X_j-x)|\cr &\leq \E
|G(X_m-X_i-x)-G(X_m-X_j-x)|\cr & \leq
\E\Big[\frac{c|X_i-X_j|^{ 2/3}}{(1+|X_m-X_i-x|^2)^{1/3}}+
   \frac{c|X_i-X_j|^{2/3}}{(1+|X_m-X_j-x|^2)^{1/3}}\Big].\cr}$$ By (5.6) and
symmetry
$$\E\Big((1+|X_m-X_i-x|^2)^{-1/2}\Big)\leq
\E\Big((1+|X_m-X_i|^2)^{-1/2}\Big)
\leq 1\land c(m-i)^{-1/2}.$$
  Then using Holder's inequality in the form $|\E( fg)|\leq \|f\|_{ 3}\,\|g\|_{
3/2}$ we obtain from the last two displays that
$$\eqalignno{\E |J(X_m-X_i-x)&-J(X_m-X_j-x)|&(9.10)\cr &
\leq
\frac{c|i-j|^{1/3}}{1\vee |m-i|^{1/3}} +\frac{c|i-j|^{1/3}}{1\vee
|m-j|^{1/3}}.\cr}$$ Thus for
$i\geq 2Bm$, summing over $j$ from $i-Bm$ to $i$ and dividing by
$Bm$ shows
$$\E|K_i-L_i|\leq c(Bm)^{1/3} (1\vee |m-i|)^{-1/3}.$$
  Therefore,
$$\eqalign{\E \sum_{i=h+2Bm}^{m} |K_i-L_i|&\leq
\sum_{i=h+2Bm}^{m}
\frac{c(Bm)^{1/3}}{(1\vee |m-i|)^{1/3}}\leq cm B^{1/3}.\cr}$$ Recalling
(9.8)--(9.9) and then using Proposition 3.1 we have that
$$\E\Big| \sum_{i=1}^{m} |K_i-L_i|\Bigr|^p\leq c (\log m)^{p}+cm^p
B^{p/3}\leq c (\log n)^{p}n^p
B^{p/3}\eqno (9.11)$$
for $n$ large. Combining with this with  (9.7), (7.5) and
Cauchy-Schwarz, the left hand side of (9.7) is bounded in $L^p$ norm by
$$c(\log n)^{k-1}  n^{1+( 1/2p)}B^{1/3}. \eqno (9.12)$$
\ms

\ms

We use summation by parts on
  $$\sum_{i=1}^{m} L_i \didi \eqno (9.13)$$  and we see that it is equal to
$$\eqalignno{&L_m [\wt B_{ k-1,n}(m,x_{ k^{ c}})-\wt B_{ k-1,n}(m,x'_{ k^{
c}})]\cr &-
\sum_{i=1}^{m} [\wt B_{ k-1,n}(i-1,x_{ k^{ c}}) -\wt B_{ k-1,n}(i-1,x'_{
k^{c}})]\, [L_i-L_{i-1}].
  &(9.14)\cr}$$  Write $w=|x_{k^c}-x_{k^c}'|/\sqrt n\leq 1$. Using the fact
that
$L_m$ is bounded by
$c\log m$ and our inductive hypothesis concerning
  (7.10), we  can bound the $L^p$ norm of the
first term of (9.14) by
$$c(\log n)^{k-1}  n^{ 1+(k-1)/p} w^{8^{ -k+1}}.$$  Since
$K_i$ is bounded by $c\log n$,  then $L_i-L_{i-1}$ is bounded by
$c\log n/(Bn)$. Hence using once again our inductive hypothesis concerning
  (7.10)
$$\eqalign{\norm{[\wt B_{ k-1,n}(i-1,x_{ k^{ c}})&-\wt B_{ k-1,n}(i-1,x'_{
k^{ c}})]\, [L_i-L_{i-1}]}_p\cr &\leq
\frac{c\log n}{Bn}  \norm{\wt B_{ k-1,n}(i-1,x_{ k^{ c}})-\wt B_{
k-1,n}(i-1,x'_{ k^{ c}})}_p\cr &\leq
\frac{c\log n}{Bn}  (\log n)^{k-2}  n^{ 1+(k-1)/p}  w^{8^{ -k+1}}.\cr}$$
Since there are
$n$ summands in the sum in (9.14), we bound the $L^p$ norm of the left
hand side of  (9.13) by
$$\frac{c}{B} (\log n)^{k-1}  n^{ 1+(k-1)/p}  w^{8^{ -k+1}}.\eqno
(9.15)$$

\ms

Notice that
$$\Bigl|\sum_{i=1}^{m} H(X_m-X_i-x) \dido\Bigr|
\leq m\wt Y_{ k-1}( n)\sup_{1\leq i\leq m} |H(X_m-X_i-x)|.\eqno
(9.16)$$ By Proposition 2.3, $H(z)$ is 0 unless $|z|\geq e^{8\log m}$. By
hypothesis we have
$|x|\leq \sqrt n$. Therefore using $|H(z)|^{ 2p}\leq c|\log z|^{ 2p}\leq
c|z|$ for
$|z|\geq e^{8\log m}$
$$\eqalign{\E |H(X_m-X_i-x)|^{ 2p}&\leq c(p) \E[ |X_m-X_i-x|;
|X_m-X_i-x|\geq  e^{8\log n}]\cr &\leq c(p)e^{-8\log m}
\E|X_m-X_i-x|^2\cr &\leq c(p) m^2 e^{-8\log m}=c(p) m^{-6}.\cr}$$
Hence $\E \sup_{ i\leq m}|H(X_m-X_i-x)|^{ 2p}\leq c(p) m^{-5}.$
Since
$w\geq 1/\sqrt n$, then this estimate, (9.16), (7.5) and Cauchy-Schwarz
imply that the left hand side of (9.16) is bounded in $L^p$ norm by
$c(\log n)^{k-2}  n^{ 1/2p}\leq c(\log n)^{k-2}  n^{1+ 1/2p}w^{2}
\leq c(\log n)^{k-1}  n^{1+ 1/2p}w^{8^{-k+1}}$.

Combining our estimates (9.12), (9.15), and our last estimate for (9.16), we
have
$$\eqalign{&\|\sum_{i=1}^{m} G(X_m-X_i-x))\didi\|_{ p}\cr &\leq c (\log
n)^{( k-1)} n^{1+( k-1)/p}
[w^{8^{ -k+1}}+B^{1/3} +w^{8^{ -k+1}} B^{-1}].
\cr}\eqno (9.17)$$
  If we take
$B=w^{6(8^{ -k})}$,  we obtain (9.6). Together with (9.2)--(9.5) we obtain
(9.1).
\qed

\subsec{10. Other results}

\ms
\ni{\sl A. $L^2$ norms.} By Section 3 of \cite{\BR} we see that we can
choose $W_t$ and $X_n$ such that $$\norm{\sup_{s\leq 1}
|X^n_s-W^n_s|\, }_2=o(n^{-\zeta})$$ for some $\zeta>0$. If we then use
this (in place of (1.7)), our proof shows that we obtain
$$\norm{ \wt \beta_k(n,0)-\wt
\gamma_k(1,0,n)}_2=o(n^{-\eta})\eqno (10.1)$$  for some $\eta>0$.

\ms
\ni{\sl B. A correction.} We take this opportunity to correct an error in
\cite{\BKb}. In the statement of (8.3) in Theorem 8.1 of that paper,
$G^\lor :=\max_{1\leq j\leq k-1} |G(x_j)|$ should be replaced by
$N^\lor := \max_{1\leq j\leq k-1} |x_j|^{-1}$. The term $G^\lor$ also needs
to be replaced by $N^\lor$ throughout the proof of (8.3).

Proposition 9.2 of that paper is correct as stated. Where the proof of this
proposition says to follow the lines of the proof of (8.3), it is to be
understood that here one uses $G^\lor$ throughout.

For the convenience of the interested reader we give a complete proof of
that proposition in the following Appendix.

\subsec{Appendix. Proof of Proposition 9.2 of \cite{\BKb}}

\ms

The proof of Proposition 9.2 in \cite{\BKb} is perhaps a bit confusing due to an
error in the statement of (8.3) in Theorem 8.1 of that paper. This Appendix
 provides a complete proof
of Proposition 9.2 of \cite{\BKb}.

\bs

Write $g(y)=\frac{1}{\pi} \log (1/|y|)$, $\wt \gamma_1(x,t)=t$, and for
$x=(x_2, \ldots, x_k)=(x_{k^c},x_k)\in (\R^2)^{k-1}$,  set
$$g^\lor(x)=\max_{2\leq i\leq k} |g(x_i)|$$ and $$\ol U_k(t,x)=
\int_0^t g(W_t-W_r-x_k) \wt \gamma_{k-1}(dr, x_{k^c}).$$

\proclaim Proposition A.1. Let $M\geq 1$, let $x,x'\in (\R^2)^{k-1}$ with
$|x_i|, |x'_i|\leq M$ for $i=2, \ldots, k$, and let
$g^+=g^\lor(x)+ g^\lor (x')+1$. There exist $a_k$  and $\nu_k$ 
 such that for $k\geq 2$
$$\E|\ol U_k(t,x)-\ol U_k(t,x')|^p\leq c (g^+)^{\nu_kp} |x-x'|^{a_kp}.
\leqno (a)$$
$$\E|\ol U_k(t,x)-\ol U_k(s,x)|^p\leq c (g^+)^{\nu_kp} |t-s|^{a_kp}.
\leqno (b)$$

\ni Except for the restriction on the size of $x, x'$, this is Proposition 9.2 of
 \cite{\BKb} translated to the notation of 
 this paper. 
Using the argument of \cite{\BKb} this is sufficient to prove the joint continuity of
$\wt \gamma_k(t,x)$ over $t\in [0,1], x\in (B(0,M))^{k-1}$ for each $k$ and $M$.
For almost every path of Brownian motion, $\{W_s: s\in [0,1]\}$ is contained
in $B(0,M)$ for some $M$ (depending on the path), and hence for $|x|>M$ we have
$\wt \gamma_k(t,x)= 0$.
The joint continuity 
of $\wt \gamma_k(t,x)$ over $t\in [0,1],
x\in (\R^2)^{k-1}$ follows.
For the purposes of  this paper we only need the case $M=1$.

Note that renormalization allows us to use $\ol U_k(t,x)$ in place
of
$$U^*_k(t,x)=\int_0^t[g(W_t-W_r-x_k)-g(-x_k)] \wt
\gamma_{k-1}(dr, x_{k^c}).$$
 If one were to try to use $U^*_k(t,x)$ in Proposition A.1, the right hand sides of
(a) and (b) would have to have $g^+$ replaced by
$N^\lor$, which is not a good enough bound for the joint continuity argument.

\proof Since $g^+$ is infinite if any component of $x$ or $x'$ is zero, we may
assume that no component of either is 0. 

 Let $A\in (0,\frac12]$ be chosen later and let
$$g_A(x)=\big(g(x)\land \tfrac{1}{\pi}\log(1/A)\big)\lor \big(-\tfrac{1}{\pi}\log(1/A)\big),$$
$$  h_A(x)=[g(x)-g_A(x)]1_{(|x|<A)}, \qq j_A(x)=[g_A(x)-g(x)]1_{(|x|>A^{-1})},$$
so $g=g_A+h_A-j_A$.
With $C_{ k}=C\cup \{ k\}$ set
$$L_{ k}(t,x)=\sum_{C\subset\{2, \ldots, k-1\}} \Big(\prod_{i\in C} |g(x_i)|\Big)
\gamma_{k-1-|C|} (t,x_{C_{ k}^c}).$$

The proof is by induction. We start with
$k=2$. In preparation for general $k$ we retain the general notation, but
note that when $k=2$,  we have
$L_{ 2}(t,x)=t$, $x_{k^c}$ is superfluous, and
we have 
$\wt \gamma_{k-1}(dr, x_{k^c})=dr$.
$$\eqalignno{\ol U_k(t,x)&-\ol U_k(t,x')&(A.1)\cr &=\int_0^t
[g_A(W_t-W_r-x_k)-g_A(W_t-W_r-x'_k)] \wt \gamma_{k-1}(dr, x_{k^c})\cr
&\qq +\int_0^t h_A(W_t-W_r-x_k) \wt \gamma_{k-1}(dr, x_{k^c})\cr &\qq 
-\int_0^t h_A(W_t-W_r-x'_k) \wt \gamma_{k-1}(dr, x_{k^c})\cr 
&\qq -\int_0^t j_A(W_t-W_r-x_k) \wt \gamma_{k-1}(dr, x_{k^c})\cr &\qq 
+\int_0^t j_A(W_t-W_r-x'_k) \wt \gamma_{k-1}(dr, x_{k^c})\cr 
&=:
I_1+I_2-I_3-I_{4}+I_{5}.\cr}$$ If we connect $x,x'$ by a curve
$\Gamma$ of length $c|x-x'|$ that never gets closer to 0 than $|x|\land |x'|$, use
the fact that $|\grad g_A|\leq A^{-1}$, and use inequality (8.1) of \cite{\BKb}
(this is only needed for $k>2$)
$$\E|I_1|^p \leq cA^{-p}|x-x'|^p \E L_k(t,x)^p \leq c A^{-p} (g^+)^{\nu_1 p} 
|x-x'|^p.$$ By Proposition 5.2 of \cite{\BKb}, for some constants
$b_1$ and $\nu'_1$
$$\E  |I_2|^p  \leq cA^{b_1p} (g^+)^{\nu'_1 p} \eqno (A.2)$$ and similarly for
$I_3$. 

We next turn to $I_{4}$. Standard estimates on Brownian motion tells
us that
$$\P(\sup_{0\leq r\leq t\leq 1} |W_t-W_r|>\lam)\leq ce^{-c'\lam^2}.
\eqno (A.3)$$
Since $|x_k|\leq M$, it follows that
$$\E[\sup_{0\leq r\leq t\leq 1} |\log(1/|W_t-W_r-x_k|)|^p]\leq c(p,M) \eqno (A.4)$$
for each $p\geq 1$. If $|A|\leq (2M)^{-1}$, then 
$|W_t-W_r-x_k|\geq A^{-1}$ only if $|W_t-W_r|\geq  (2A)^{-1}$.
So by Cauchy-Schwarz, (A.3), and (A.4),
$$\eqalign{\E[&\sup_{0\leq r\leq t\leq 1} |j_A(W_t-W_r-x_k)|^p\cr
&\leq c\E[\sup_{0\leq r\leq t\leq 1} |\log (1/|W_t-W_r-x_k|)|^{p}
1_{(\sup_{0\leq r\leq t\leq 1} |W_t-W_r-x_k|\geq A^{-1})}]\cr
&\leq c(\E[\sup_{0\leq r\leq t\leq 1} |\log (1/|W_t-W_r-x_k|)|^{2p}])^{1/2}
(\P(\sup_{0\leq r\leq t\leq 1} |W_t-W_r|\geq  (2A)^{-1}))^{1/2}\cr
&\leq cA^p\cr}$$
for each $p\geq 1$.
Using the fact that 
$$|I_{4}|\leq \big(\sup_{0\leq r\leq t\leq 1} |j_A(W_t-W_r-x_k)|\big) L_k(t,x),
$$
another application of Cauchy-Schwarz shows that
$$\E |I_{4}|^p \leq (g^+)^{\ol \nu_1 p} A^p.$$
$I_{5}$ is handled the same way.

Combining shows the left hand side of (A.1) is bounded by 
$$ c(g^+)^{\nu_2 p}[A^{-p} |x-x'|^{p}+ A^{b_1p}+A^p]$$ for some constant
$\nu_2$, and we obtain (a) for $k=2$ by setting $A=|x-x'|^{1/2} \land
 (2M)^{-1}$.

Next we look at (b) for the $k=2$ case. We write
$$\eqalignno{\ol U_k(t,x)&-\ol U_k(s,x)&(A.5)\cr &=\int_0^s [g_A(W_t-W_r-x_k)
- g_A(W_s-W_r-x_k)] \wt \gamma_{k-1}(dr, x_{k^c})\cr &\qq+\int_s^t
g_A(W_t-W_r-x_k) \wt \gamma_{k-1}(dr, x_{k^c})\cr &\qq +\int_0^t
h_A(W_t-W_r-x_k) \wt \gamma_{k-1}(dr, x_{k^c})\cr 
&\qq-\int_0^s h_A(W_s-W_r-x_k) \wt \gamma_{k-1}(dr, x_{k^c})\cr 
&\qq-\int_0^s j_A(W_s-W_r-x_k) \wt \gamma_{k-1}(dr, x_{k^c})\cr 
&\qq+\int_0^s j_A(W_s-W_r-x'_k) \wt \gamma_{k-1}(dr, x_{k^c})\cr 
&=:
I_6+I_7+I_8-I_9-I_{10}+I_{11}.\cr}$$ Using  Cauchy-Schwarz, for
$s,t\leq 1$ we have for some constant $\nu''_1$
$$\E |I_6|^p \leq cA^{-p} s^p \Big(\E |W_t-W_s|^{2p}\Big)^{1/2}
\Big(\E L(s,x_{k^c})^{2p}\Big)^{1/2}\leq cA^{-p}(g^+)^{\nu''_1 p}|t-s|^{p/2}.$$
We bound $I_7$ by
$$\E|I_7|^p\leq c (\log(1/A))^p |t-s|^p (g^+)^{\nu''' p}.$$ We bound
$I_8$ and $I_9$ just as we did $I_2$ and bound $I_{10}$ and $I_{11}$ as we did 
$I_{4}$.  Combining, the left hand side of (A.5) is
bounded by
$$c(g^+)^ {\nu_2 p} [A^{-p}|t-s|^{p/2} + A^{b_1p}+A^p],$$ and (b) follows by setting
$A=|t-s|^{1/4}\land (2M)^{-1}$.

We now turn to the case when $k>2$. We suppose (a) and (b) hold for
$k-1$ and prove them for $k$.  We prove (a) in two cases, when
$x_{k^c}=x'_{k^c}$ and when $x_k=x'_k$; the general case follows by the
triangle inequality. Suppose first that $x_{k^c}=x_{k^c}'$. Using the induction
hypothesis, the proof is almost exactly the same as the proof of (a) in the case
$k=2$. 

Suppose next that $x_k=x'_k$. Let 
$$V_A=\{|W_{s+u}-W_u|\geq u^{1/4}/A \hbox{ for some } s,u\in [0,1]\}.$$
Standard estimates on Brownian motion show that
$$\P(V_A)\leq c_1 e^{-c_2/A^2}.$$ We write
$$\eqalignno{\ol U_k(t,x)&-\ol U_k(t,x')&(A.6)\cr &=1_{V_A}
\int_0^t g_A(W_t-W_r-x_k) \wt \gamma_{k-1}(dr, x_{k^c})\cr &\qq-1_{V_A}
\int_0^t g_A(W_t-W_r-x_k) \wt \gamma_{k-1}(dr, x'_{k^c})\cr 
&\qq+  \int_0^t h_A(W_t-W_r-x_k)\wt \gamma_{k-1}(dr, x_{k^c})\cr 
&\qq-  \int_0^t h_A(W_t-W_r-x_k)\wt \gamma_{k-1}(dr, x'_{k^c})\cr 
&\qq-  \int_0^t j_A(W_t-W_r-x_k)\wt \gamma_{k-1}(dr, x'_{k^c})\cr 
&\qq+  \int_0^t j_A(W_t-W_r-x'_k)\wt \gamma_{k-1}(dr, x'_{k^c})\cr 
&\qq+1_{V_A^c}
\int_0^t g_A(W_t-W_r-x_k) [\wt
\gamma_{k-1}(dr, x_{k^c})-\wt \gamma_{k-1}(dr, x'_{k^c})]\cr &=:
I_{12}-I_{13}+I_{14}-I_{15}-I_{16}+I_{17}+I_{18}.\cr}$$ Since $|g_A|\leq \log (1/A)$, for some
constant $\nu_{k-1}$
$$\eqalign{\E|I_{12}|^p&\leq c(\log (1/A))^p \E[L_k(t, x)^p;V_A]\cr &\leq c(\log
(1/A))^p \Big(\E L_k(t, x)^{2p}\Big)^{1/2}\Big(
\P(V_A)\Big)^{1/2}
\leq c(g^+)^{\nu_{k-1}p} A^p,\cr}$$ and similarly for $I_{13}$. We bound
$I_{14}$ and $I_{15}$ just as we did $I_2$ and bound $I_{16}$ and $I_{17}$ as
we did $I_{4}$.

We turn to $I_{18}$. Let $f(r)=g_A(W_t-W_r-x_k)$ and 
$$f_A(t)=\frac{1}{A^{12}} \int_t^{t+A^{12}} f(u) \, du.$$ On $V_A^c$ we have 
$$\eqalign{|f(r)-f(s)|&=|g_A(W_t-W_r-x_k)-g_A(W_t-W_s-x_k)|\cr &\leq A^{-1}
|W_r-W_s|\leq A^{-2} |r-s|^{1/4},\cr}$$ and therefore
$$|f(t)-f_A(t)|\leq A^{-2} (A^{12})^{1/4} =A.$$ Using integration by parts, we
write
$$\eqalign{I_{18}&=1_{V^c_A}\int_0^t [f(r) -f_A(r)] \wt
\gamma_{k-1}(dr, x_{k^c})\cr &\qq - 1_{V^c_A}\int_0^t [f(r) -f_A(r)]
\wt \gamma_{k-1}(dr, x'_{k^c})\cr &\qq +  1_{V_A^c}f_A(t) [\wt
\gamma_{k-1}(t, x_{k^c})-\wt \gamma_{k-1}(t, x'_{k^c})]\cr &\qq
-1_{V_A^c}\int_0^t 
   [\wt \gamma_{k-1}(r, x_{k^c})-\wt \gamma_{k-1}(r, x'_{k^c})] f_A(dr)\cr &=:
I_{19}-I_{20} +I_{21}-I_{22}.\cr}$$ We bound
$$\E |I_{19}|^p\leq cA^p (g^+)^{\nu'_{k-1} p}$$ for some constant
$\nu'_{k-1}$ and similarly for $I_{20}$. By the induction hypothesis and the
fact that
$|f_A|$ is bounded by $\log (1/A)\leq cA^{-p}$ ,
$$\E |I_{21}|^p\leq c A^{-p} (g^+)^{\nu''_{k-1} p} |x_{k^c}-x'_{k^c}|^{a_{k-1}
p}.$$

Finally, since
$|f'_A| \leq \norm{f_A}_\infty A^{-12}\leq cA^{-13}$,
$$\E |I_{22}|^p\leq c(g^+)^{\nu'''_{k-1} p} A^{-13p}
|x_{k^c}-x'_{k^c}|^{a'_{k-1}p}.$$ If we combine all the terms, we see that the left
hand side of (A.6) is bounded by
$$c(g^+)^{\nu_k p}[ A^{b_kp} + A^{-b'_kp} |x-x'|^{a''_{k-1}p}+A^p].$$ Setting $A=
|x-x'|^{a''_{k-1}/(2b'_k)}\land (2M)^{-1}$ completes the proof of (b) for $k>2$.

The proof of (b) for $k>2$ is almost identical to the  $k=2$ case.
\qed

\bs

\bs

\subsec{References}

\refer{\Bas} R.F. Bass, {\sl Probabilistic Techniques in Analysis}. Springer,
New York, 1995.

\refer{\BKa} R.F. Bass and D. Khoshnevisan, Local times on curves and
uniform invariance principles. {\sl Prob. Th. and rel. Fields \bf  92} (1992)
465--492.

\refer{\BKb} R.F. Bass and D. Khoshnevisan, Intersection local times and
Tanaka formulas. {\sl  Ann.  Inst. H. Poincar\'e \bf 29} (1993) 419--451.

\refer{\BKc} R.F. Bass and D. Khoshnevisan, Strong approximations to
Brownian local time, {\sl Seminar on Stochastic Processes,  1992}, 43--66,
Birkh\"auser, Boston, 1993.

\refer{\BR} R.F. Bass and J.~Rosen,  An almost sure invariance principle
   for the range of  planar random walks, preprint.

\refer{\BP} I. Berkes and  W.  Philipp,  Approximation theorems for
independent and weakly dependent random vectors,  {\sl Ann. Probab.
\bf 7} (1979) 29--54.

\refer{\Da} E. B.~Dynkin,  Self-intersection gauge for random walks and for
Brownian motion. {\sl Ann. Probab.
\bf 16} (1988) 1--57.

\refer{\Db} E. B.~Dynkin,  Regularized self-intersection local times of planar
Brownian motion. {\sl Ann. Probab. \bf 16} (1988) 58--74.

\refer{\LGa} J.-F.~Le Gall,  Wiener sausage and self intersection local times.
{\sl J. Funct.\ Anal. \bf  88} (1990) 299--341.

\refer{\LGb}
        J.-F. {Le Gall}.
       Propri\'et\'es d'intersection des marches al\'eatoires, I.
         Convergence vers le temps local d'intersection,
    {\sl Comm.\ Math.\ Phys.  \bf  104} (1986) 471--507.

\refer{\LGc}
      J.-F. {Le Gall}.  Some properties of planar {B}rownian motion, {\sl
\'Ecole d'\'et\'e de
      probabilit\'es de St. Flour XX, 1990}
        Springer-Verlag, Berlin, 1992.

  \refer{\La} G. Lawler, {\sl Intersections of random walks}. {Birkh\"auser},
Boston, 1991.

\refer{\Lw} G. Lawler, in preparation.

\refer{\MR} M.~Marcus and J.~Rosen, {\sl Renormalized self-intersection
local times and
       {W}ick power chaos processes}, Memoirs of the AMS, Volume 142,
Number 675, Providence, AMS, 1999.

\refer{\Pt} V.V. Petrov,  {\sl Sums of Independent Random Variables}.
Springer-Verlag, Berlin, 1975.

\refer{\Ra} J.~Rosen,  Random walks and intersection local time.
       {\sl Ann. Probab. \bf  18} (1990) 959--977.

\refer{\Rb} J.~Rosen,  Joint continuity of renormalized intersection local
times, {\sl Ann.\
        Inst.\ H.~Poincar{\'e} \bf  32} (1996) 671--700.

\refer{\Rc}
       J.~Rosen. Dirichlet processes and an intrinsic characterization for
renormalized intersection local times, {\sl Ann.\
        Inst.\ H.~Poincar{\'e} \bf 37} (2001) 403--420.

\refer{\Sp} F. Spitzer, {\sl Principles of Random Walk, 2nd ed.}, Springer,
New York, 1976.

\refer{\Va}
       S.~R.~S. Varadhan.
       Appendix to {E}uclidean quantum field theory by {K}. {S}y\-man\-zyk.
        In R.~Jost, editor,  {\sl Local Quantum Theory}. Academic Press,
Reading, MA,
      1969.

\bs
\ni R. Bass Address: Department of Mathematics, University of Connecticut,
\hfill\break Storrs, CT 06269-3009 (bass@math.uconn.edu)
\ms
\ni J. Rosen Address: Department of Mathematics, College of Staten Island,
CUNY,
\hfill\break Staten Island, NY, 10314 (jrosen3@earthlink.net)

\bye